\documentclass[a4paper,12pt]{amsart}
\usepackage{amsfonts}
\def\T{{\mathcal{T}}}
\def\M{{\mathcal{M}}}
\def\O{{\mathcal{O}}}
\def\V{{\mathcal{V}}}
\def\N{{\mathbb{N}}}
\def\oM{{{\overline{\mathcal{M}}}}}
\def\oH{{\overline{H}}}
\def\hH{{\widehat{H}}}
\def\wS{{\widehat{S}}}
\def\oll{{\overline{ll}}}
\def\hll{{\widehat{ll}}}
\def\CP1{{\mathbb{C}\mathrm{P}^1}}

\def\E{{\mathbb{E}}}
\def\L{{\mathbb{L}}}
\def\vphi{\varphi}
\def\aut{\mathrm{aut}}
\def\l{\langle}
\def\r{\rangle}
\def\ll{\langle\langle}
\def\rr{\rangle\rangle}
\def\d{\partial}
\def\res{\mathrm{res}}

\newtheorem{theorem}{Theorem}
\newtheorem{lemma}{Lemma}

\title[Meromorphic functions and moduli spaces of curves]{Geometry of meromorphic functions and intersections on moduli spaces of curves}

\author{S.~Shadrin}

\begin{document}

\begin{abstract}
In this paper we study relations between intersection numbers on moduli spaces
of curves and Hurwitz numbers. First, we prove two formulas expressing
Hurwitz numbers of (generalized) polynomials via intersections on moduli spaces
of curves. Then we show, how intersection numbers can be expressed via Hurwitz
numbers. And then we obtain an algorithm expressing intersection numbers
$\l\tau_{n,m}\prod_{i=1}^{r-1}\tau_{0,i}^{k_i}\r_g$ via correlation functions of primaries.
%
%
\end{abstract}

\maketitle

\tableofcontents




\section{Introduction}

In~\cite{i}, Ionel developed a very beautiful approach to study intersection theory of moduli space of curves. Roughly speaking, the situation is the following. Consider the space of meromorphic functions with fixed genus, degree, and ramification type. There are two mappings of this space. One mapping ($ll$) takes a meromorphic function to its target curve (of genus zero) with marked critical values. Another mapping ($st$) takes a meromorphic function to its domain curve (of genus $g$) with marked critical points. Then one can relate intersection theories in the images of these mappings. This idea was also used in~\cite{p,bp} in low genera.

In this paper, we just study several applications of Ionel's technique. Concretely, we express Hurwitz numbers via intersection numbers and vice versa. More or less, the same problem, but in much more general case, is studied by Okounkov and Pandharipande in~\cite{op,op2}.

Our results could be split into three parts, which we will describe now.

\subsection{Hurwitz numbers of polynomials}
Consider polynomials with fixed critical values and fixed ramification type over each critical value. Roughly speaking, Hurwitz number is the number of such nonequivalent polynomials, or it is better to say that Hurwitz number is the multiplicity of the corresponding mapping $ll$. There are several combinatorial formulas for such numbers, see~\cite{gj,lz,l,pz,z}. In this paper, we express these numbers in terms of intersection numbers on the moduli space of curves of genus zero with marked points.

Ionel's theory (Sections~\ref{admissible-covers} and~\ref{lemma-of-ionel}) reduces the problem of counting a Hurwitz number to the problem of calculating the homology class of the image of the mapping $st$. We do this for polynomials in Section~\ref{poincare-dual}.

Really, the formula we have obtained is rather complicated (Theorem~\ref{thm1}. It uses some recursively defined classes (Section~\ref{coh-classes}) and it is hard to work with even in the simplest cases (see examples in Section~\ref{independent-check}).

Nevertheless, we think that this formula is beautiful itself, and some steps of its proof are excellent examples of the technique we use in this paper.

\subsection{Hurwitz numbers of generalized polynomials and two-pointed
ramification cycles}

Consider a space of meromorphic functions defined on genus $g$ curves with a fixed collection of ramification data. If only two critical values of these functions are not simple, then the homology class of the image of this space under the mapping $st$ is called two-pointed ramification cycle.

The situation with Hurwitz numbers (multiplicities of the mapping $ll$) in this case is just the same as for polynomials; we just have to calculate the two-pointed ramification cycle. But in the case of an arbitrary $g$, it is a hard problem.

We consider the case, when one nonsimple critical value is arbitrary and the other one is the value at a point of total ramification. Such functions defined on curves of arbitrary genus we call generalized polynomials. In this particular case we did not manage to compute the corresponding two-pointed ramification cycle, but we found a way to simplify it.

Thus we obtain a formula expressing Hurwitz numbers of generalized polynomials via simplified two-pointed ramification cycles (Theorem~\ref{theorem-generalized-polynomials}). This formula also could appear to be useless, but we found a remarkable application of it.

Using the technique of Ionel many times (Lemmas~\ref{S-expression-via-difference} and~\ref{lemma-l-independence}), we found an expression for some intersection numbers via these simplified two-pointed ramification cycles (Lemma~\ref{lemma-proof2}). This gives us an expression of these intersection numbers via Hurwitz numbers of generalized polynomials (Theorem~\ref{formulas-for-tau3g}).

We consider the intersection numbers which look as follows. Let $\oM_{g,1}$ be the compactified moduli space of genus $g$ curves with one marked point. By $\psi$ denote the first Chern class of the line bundle over $\oM_{g,1}$, whose fiber at a point of $\oM_{g,1}$ is the cotangent line at the marked point of the corresponding curve. We consider the intersection numbers
\begin{equation}
\l\tau_{3g-2}\r_g=\int_{\oM_{g,1}}\psi^{3g-2}.
\end{equation}

Previously we knew the unique way to calculate such numbers. By the Witten--Kontsevich theorem these numbers are equal to the coefficients of the string solution of the KdV hierarchy, and it is easy to compute these coefficients. Our formula leads to a combinatorial way to calculate such
numbers. We express these intersection numbers via Hurwitz numbers, and it is a purely combinatorial problem to calculate a Hurwitz number.

In addition our formula gives an infinite number of linear relations for Hurwitz number, which can also be useful.

\subsection{A way to compute
$\l\tau_{n,m}\prod_{i=1}^{r-1}\tau_{0,i}^{k_i}\r_g$}

The intersection numbers $\l\tau_{n,m}\prod_{i=1}^{r-1}\tau_{0,i}^{k_i}\r_g$ are the natural generalization of the numbers $\l\tau_{3g-2}\r_g$. The definition of these numbers is rather complicated, so we do not recall it in the inroduction.

For these numbers, we try to do just the same as for $\l\tau_{3g-2}\r_g$. First, we generalize Lemma~\ref{lemma-proof2} and we obtain an expression for these numbers in terms of the integrals against two-pointed ramification cycles (Theorem~\ref{algorithm1}). This time these integrals do not correspond to any Hurwitz numbers, but they are very similar to Hurwitz numbers. In particular, we found a generalization of a standard recursive relation for Hurwitz numbers, which works for our integrals (Theorem~\ref{algorithm2}).

Theorem~\ref{algorithm2} reduces calculation of any of our integrals against two-pointed ramification cycles to calculation of simple intersection numbers in genera zero and one. The
computing of these simple intersection numbers in genus one is not completely clear, but we discuss how it can be reduced to the intersection numbers in genus zero (Section~\ref{genusone}).

The motivation to study the intersection numbers $\l\tau_{n,m}\prod_{i=1}^{r-1}\tau_{0,i}^{k_i}\r_g$ looks as follows. These numbers conjectured by E.~Witten to coinside with some coefficients of the string solution of the $r$th Gelfand--Dikii hierarchy (we recall this conjecture in Section~11). So, using our theorems, one can compute these numbers and compare them with the coefficients of the string solution. Thus the Witten's conjecture could be checked in particular cases.

In the appendix, we compute the intersection number $\l\tau_{6,1}\r_3$ in the case $r=3$ in two ways, and thus we check the Witten's conjecture in a very particular case.

\subsection{Organization of the paper}

In Section~2, we define Hurwitz numbers. In Section~3, we formulate our first theorem expressing Hurwitz numbers of usual polynomials with arbitrary ramification via intersections on $\oM_{0,n}$.
In Section~4, we recall the definition of admissible covers. In Section~5, we formulate the lemma of E.~Ionel which plays the principal role in our paper. In Section~6, we prove our first theorem.

In Sections~7 and~8, we formulate and prove our second theorem. There we give an
expression for Hurwitz numbers of generalized polynomials with one nonsimple
critical value. In Sections~9 and~10, we express the intersection number
$\l\tau_{3g}\tau_{0}^2\r_g$ via Hurwitz numbers.

In Section~11, we recall the definition of the intersection numbers $\l\tau_{n,m}\prod_{i=1}^{r-1}\tau_{0,i}^{k_i}\r_g$. There we also recall the conjecture of E.~Witten. In Sections~12 and~13, we present an algorithm for calculating all intersection numbers of the type $\l\tau_{n,m}\prod_{i=1}^{r-1}\tau_{0,i}^{k_i}\r_g$. In the appendix, we show an example of usage of this algorithm.

\subsection{Acknowledgments} The author is grateful to S.~K.~Lando,
S.~M. Natanzon, and M.~Z.~Shapiro for useful remarks and discussions.




\section{Definition of Hurwitz numbers}

In this section, we give a definition of Hurwitz numbers (Section~\ref{def-Hur-num}). To define a Hurwitz number, we must fix a collection of passports whose definition we give in Section~\ref{def-pass}.

\subsection{Definition of passports}\label{def-pass}

Consider a meromorphic function $f\colon C\to\CP1$ of degree $n$ defined on a smooth curve $C$ of genus $g$. Let $z$ be a point of $\CP1$. Then $f^{-1}(z)=a_1p_1+\dots+a_lp_l$, where $p_1,\dots,p_l$ are pairwise distinct points of $C$ and $a_1,\dots,a_l$ are positive integers such that $\sum_{i=1}^la_i=n$. Suppose that $a_1\geq a_2\geq\dots\geq a_l$. Then the tuple of numbers $(a_1,\dots,a_l)$ is called the \emph{passport of $f$ over $z$}.

For instance, the passport of $f$ over a regular point is equal to $(1,\dots,1)$ and the passport of $f$ over a simple critical value is equal to $(2,1,\dots,1)$. If $f$ is a polynomial of degree $n$, then the passport of $f$ over infinity is equal to $(n)$.

\subsection{Definition of Hurwitz numbers}\label{def-Hur-num}

Two meromorphic functions, $f_1\colon C_1\to\CP1$ and $f_2\colon C_2\to\CP1$, are isomorphic if there exists a biholomorphic map $\vphi\colon C_1\to C_2$ such that $f_1=f_2\vphi$.

Consider $m$ distinct points of $\CP1$. Let $A_i=(a^i_1,\dots,a^i_{l_i})$, $i=1,\dots,m$, be nonincreasing sequences of positive integers such that $\sum_{j=1}^{l_i}a^i_j=n$. Up to isomorphism, there is a finite number of meromorphic functions $f\colon C\to\CP1$ of degree $n$ defined on smooth curves of genus $g$ such that the passport of $f$ over $z_i$ is equal to $A_i$, $i=1,\dots,m,$ and $f$ is unramified over $\CP1\setminus\{z_1,\dots,z_m\}$.

The Hurwitz number $h(g,n|A_1,\dots,A_m)$ is the weighted count of such functions, where a function $f\colon C\to\CP1$ is weighted by $1/|\aut(f)|$. For example, the number $h(g,2|(2),\dots,(2))$ is equal to $1/2$.

In~\cite{a, elsv, gv, l}, one can find some examples of formulas for Hurwitz numbers in special cases.




\section{Hurwitz numbers of polynomials}

The goal of this section is to state our formula for Hurwitz numbers of polynomials (Section~\ref{formula-hur-pol}). In Section~\ref{intro}, we fix notations for the ramification data and recall the known combinatorial formula for Hurwitz numbers of polynomials. In Section~\ref{coh-classes}, we define the cohomological classes which we use in our formula.

\subsection{Introduction}\label{intro}

Let us fix $n\in\N$, $g=0$, and a collection of passports $A_1,\dots,A_m$ with usual requirements (for any $i$ we have $A_i=(a^i_1,\dots,a^i_{l_i})$, $a^i_1\geq\dots\geq a^i_{l_i}$, and
$\sum_{j=1}^{l_i}a^i_j=n$). We assume that $l_1=1$ and $A_1=(n)$. Also we assume that $\sum_{i=1}^m\sum_{j=1}^{l_i}(a^i_j-1)=2n-2$ (the Riemann--Hurwitz formula).

The following formula is proved in~\cite{gj, z, lz}:
\begin{equation}\label{formula-gjlz}
h(0,n|A_1,\dots,A_m)=n^{m-3}
\cdot \frac
{|\aut(l_2,\dots,l_m)|}
{|\aut(A_2,\dots,A_m)|} \cdot \prod_{i=2}^m \frac
{(l_i-1)!} {|\aut(A_i)|}.
\end{equation}

In this paper, we prove another formula expressing the same numbers via intersections.

\subsection{Cohomological classes on $\oM_{0,N}$}\label{coh-classes}

In Section~\ref{standard}, we fix notations and give the standard definition of $\psi$-classes on the moduli space of curves. In Section~\ref{psi-p-classes}, we give a rather complicated definition of very specific $\Psi_p$-classes. We need this just to write down our formula for Hurwitz numbers of polynomials in a compact way. In Section~\ref{well-defined}, we prove that our
$\Psi_p$-classes are well-defined. This is necessary since the definition of $\Psi_p$-classes is a recursive one and it is not obvious that the recursive relations are compatible.

\subsubsection{Standard definitions}\label{standard}

Let $N$ be equal to $\sum_{j=1}^m l_j$. We consider the moduli space of genus zero curves with $N$ marked points $\oM_{0,N}\ni (C, x_1^1, x^2_1,\dots,x^2_{l_2}, \dots, x^m_1,\dots,x^m_{l_m})$; we mean a one-to-one correspondence between marked points and indices of $a^i_j$.

By $\pi_{p,q,k}$ we denote the projection $\oM_{0,N}\to\oM_{0,2+l_p}$ that forgets all points except for $x_1^1$, $x^p_1,\dots,x^p_{l_p}$, and $x_k^q$.

Consider the moduli space of genus zero curves with $k$ marked points $\oM_{0,k}\ni (C, y_1, \dots, y_k)$. By $\psi(y_j)$ we denote the first Chern class of the line bundle over $\oM_{0,k}$ whose fiber at the point $(C, y_1, \dots, y_k)$ is the cotangent line $T^*_{y_j}C$.

For example, $\psi(x^p_i)$ stands for the first Chern class of the corresponding cotangent line bundle over $\oM_{0,N}\ni (C, x_1^1,x^2_1,\dots,x^m_{l_m})$ as well as for the first Chern class of the corresponding cotangent line bundle over $\oM_{0,2+l_p}\ni (C, x_1^1, x^p_1,\dots,x^p_{l_p}, x_k^q)$.

\subsubsection{Definition of $\Psi_p$-classes}\label{psi-p-classes}

We define classes $\Psi_p(b^i_j)^{i\in\{2,\dots,m\}\setminus\{p\}}_{j=1,\dots,l_i}$ depending on
$N-l_p-1$ indices corresponding to all points $x^i_j$ except for $x^1_1,x^p_1,\dots,x^p_{l_p}$. If $b^i_j=0$ for all $i$ and $j$, then we put $\Psi_p(b^i_j)=0$.

We further give the following recursive definition. Suppose we have already
defined all classes with
\begin{equation}
\sum_{i\in\{2,\dots,m\}\setminus\{p\}}\sum_{j=1}^{l_i}b^i_j\leq s.
\end{equation}
Consider a sequence $(b^i_j)^{i\in\{2,\dots,m\}\setminus\{p\}}_{j=1,\dots,l_i}$ such that
\begin{equation}
\sum_{i\in\{2,\dots,m\}\setminus\{p\}}\sum_{j=1}^{l_i}b^i_j=s.
\end{equation}
We fix $q\in\{2,\dots,m\}\setminus\{p\}$ and $k\in\{1,\dots,l_q\}$. Then we define $(\widehat b^i_j)^{i\in\{2,\dots,m\}\setminus\{p\}}_{j=1,\dots,l_i}$ in the following way. We put $\widehat b^q_k=\widehat b^q_k+1$; for all other indices $i,j$ we put $\widehat b^i_j=b^i_j$.

The formula for $\Psi_p(\widehat b^i_j)$ is the following one:
\begin{equation}\label{definition-of-psi-classes}
\Psi_p(\widehat b^i_j)=\widehat b^q_k \pi_{p,q,k}^*\psi(x^q_k)\Psi_p(b^i_j)-
\sum_{U} a_UD_{U\cup\{x^q_k\}}\Psi_p((b_U)^i_j)
\end{equation}
Here $a_U=\sum_{x^i_j\in U}b^i_j$; if $x^i_j\not\in U\cup\{x^q_k\}$, then $(b_U)^i_j=b^i_j$; if $x^i_j\in U$, then $(b_U)^i_j=0$; and $(b_U)^q_k=a_U+b^q_k$. We take the sum over all $U\subset\{x^i_j\}\setminus\{x^1_1,x^p_1,\dots,x^p_{l_p},x^q_k\}$.

By $D_V$, where $V\subset\{x^i_j\}$, we denote the Poincar\'e dual of the cycle, defined by the divisor in $\oM_{0,N}$ whose generic point is represented by a two-component curve such that all points from $U$ lie on one component and all points from $\{x^i_j\}\setminus U$ lie on the other component.

In the foregoing formulas, we will use only classes $\Psi_p=\Psi_p(b^i_j)$, where $b^i_j=a^i_j-1$.

\subsubsection{The $\Psi_p$-classes are well-defined}\label{well-defined}

We have to prove that our definition of $\Psi_p(b^i_j)$-classes is correct.

We fix $p$, $x^{q_1}_{k_1}$, $x^{q_2}_{k_2}$, and $(b^i_j)^{i\in\{2,\dots,m\}\setminus\{p\}}_{j=1,\dots,l_i}$. Let $\widehat b^{q_1}_{k_1}=b^{q_1}_{k_1}+1$, $\widehat b^{q_2}_{k_2}=b^{q_2}_{k_2}+1$, and $\widehat b^i_j=b^i_j$ for all other $i$ and $j$.

If we apply formula~(\ref{definition-of-psi-classes}) twice, first for $x^{q_1}_{k_1}$ and then for $x^{q_2}_{k_2}$, we obtain the following expression:
\begin{multline}\label{double-expression-for-psi-classes}
\Psi_p(\widehat b^i_j)=\widehat b^{q_1}_{k_1}\widehat
b^{q_2}_{k_2}\pi_{p,q_1,k_1}^*\psi(x^{q_1}_{k_1})
\pi_{p,q_2,k_2}^*\psi(x^{q_2}_{k_2})\Psi_p(b^i_j)\\
-\sum_{U}(\sum_{x^i_j\in U}b^i_j)\widehat
b^{q_2}_{k_2}D_{U\cup\{x^{q_1}_{k_1}\}}
\pi_{p,q_2,k_2}^*\psi(x^{q_2}_{k_2})
\Psi_p(\dots)\\
-\sum_{U}(b^{q_2}_{k_2}+\sum_{x^i_j\in U}b^i_j)\widehat
b^{q_2}_{k_2}D_{U\cup\{x^{q_1}_{k_1},x^{q_2}_{k_2}\}}
\pi_{p,q_2,k_2}^*\psi(x^{q_2}_{k_2})
\Psi_p(\dots)\\
-\sum_U(\sum_{x^i_j\in U}b^i_j)\widehat
b^{q_1}_{k_1}D_{U\cup\{x^{q_2}_{k_2}\}}
\pi_{p,q_1,k_1}^*\psi(x^{q_1}_{k_1})
\Psi_p(\dots)\\
-\sum_{U}(\widehat b^{q_1}_{k_1}+\sum_{x^i_j\in U}b^i_j)
(b^{q_1}_{k_1}+\widehat b^{q_2}_{k_2}+\sum_{x^i_j\in U}b^i_j)
D_{U\cup\{x^{q_1}_{k_1},x^{q_2}_{k_2}\}}
\pi_{p,q_1,k_1}^*\psi(x^{q_1}_{k_1})
\Psi_p(\dots)\\
+\sum_{U}(\sum_{x^i_j\in U}b^i_j)D_{U\cup\{x^{q_2}_{k_2}\}}
\sum_{V}(\sum_{x^i_j\in V}b^i_j)D_{V\cup\{x^{q_1}_{k_1}\}}
\Psi_p(\dots)\\
+\sum_{U}(\sum_{x^i_j\in U}b^i_j)D_{U\cup\{x^{q_2}_{k_2}\}}
\sum_{V}(b^{q_2}_{k_2}+\sum_{x^i_j\in U\cup
V}b^i_j)D_{U\cup V\cup\{x^{q_1}_{k_1},x^{q_2}_{k_2}\}} \Psi_p(\dots)\\
+\sum_{U}(\widehat b^{q_1}_{k_1}+\sum_{x^i_j\in U}b^i_j)
D_{U\cup\{x^{q_1}_{k_1},x^{q_2}_{k_2}\}}
\sum_{V}(\sum_{x^i_j\in V}b^i_j)
D_{U\cup V\cup\{x^{q_1}_{k_1},x^{q_2}_{k_2}\}}
\Psi_p(\dots)
\end{multline}
Here we take sums over all $U\subset\{x^i_j\} \setminus\{x^1_1,x^p_1,\dots,x^p_{l_p},x^{q_1}_{k_1},x^{q_2}_{k_2}\}$
and over all $V\subset\{x^i_j\} \setminus\{x^1_1,x^p_1,\dots,x^p_{l_p},x^{q_1}_{k_1},x^{q_2}_{k_2}\}$
such that $U\cap V=\emptyset$. We do not write down the indices in $\Psi_p$-classes since it is very easy to reconstruct these indices from the rest of each summand.

Now, to ensure that the $\Psi_p$-classes are well-defined, it is enough to prove the following lemma.

\begin{lemma}
The right-hand side of formula~(\ref{double-expression-for-psi-classes}) is
symmetric with respect to the changes $x^{q_1}_{k_1}\leftrightarrow x^{q_2}_{k_2}$ and
$b^{q_1}_{k_1}\leftrightarrow b^{q_2}_{k_2}$.
\end{lemma}

\begin{proof} Obviously, this lemma can be reduced to the fact that the sum of the last two terms of the expression is symmetric. We prove it as follows.

We fix the partition $U\cup V$. On the divisor $D=D_{U\cup V\cup \{x^{q_1}_{k_1},x^{q_2}_{k_2}\}}$, the first Chern class of the cotangent line at the double point is equal to
\begin{equation}
\psi(*)D=D\cdot\sum_{W\ni x^\alpha_\beta}(D_{W\cup\{x^{q_2}_{k_2}\}}+
D_{W\cup\{x^{q_1}_{k_1},x^{q_2}_{k_2}\}})
=\sum_WD_{W\cup\{x^{q_1}_{k_1}
x^{q_2}_{k_2}\}}D.
\end{equation}
Of course $\psi(*)$ is symmetric.

Note that
\begin{multline}
(b^{q_2}_{k_2}+\sum_{x^i_j\in U\cup
V}b^i_j)(\sum_{x^i_j\in U\cup
V}b^i_j)\psi(*)D=\\
(b^{q_2}_{k_2}+\sum_{x^i_j\in U\cup
V}b^i_j)\sum_{x^\alpha_\beta\in U\cup V}b^{\alpha}_\beta\psi(*)D,
\end{multline}
\begin{multline}
(b^{q_2}_{k_2}+\sum_{x^i_j\in U\cup
V}b^i_j)\sum_{x^\alpha_\beta\in U\cup V}b^{\alpha}_\beta\psi(*)D=\\
\sum_{W}(\sum_{x^i_j\in W}b^i_j)
(b^{q_2}_{k_2}+\sum_{x^i_j\in U\cup V}b^i_j)D_{W\cup\{x^{q_2}_{k_2}\}}D\\
+\sum_{W}(\sum_{x^i_j\in W}b^i_j)
(b^{q_2}_{k_2}+\sum_{x^i_j\in U\cup V}b^i_j)D_{W\cup\{x^{q_1}_{k_1}
x^{q_2}_{k_2}\}}D,
\end{multline}
\begin{equation}
(b^{q_1}_{k_1}+1)(\sum_{x^i_j\in U\cup V}b^i_j)\psi(*)D=
(b^{q_1}_{k_1}+1)(\sum_{x^i_j\in U\cup V}b^i_j)\sum_WD_{W\cup\{x^{q_1}_{k_1}
x^{q_2}_{k_2}\}}D.
\end{equation}

Thus we see that the sum of the last two terms of expression~(\ref{double-expression-for-psi-classes}) (with fixed $U\cup V$) is
equal to
\begin{multline}
\left((b^{q_1}_{k_1}+b^{q_2}_{k_2}+1)(\sum_{x^i_j\in U\cup V}b^i_j)+
(\sum_{x^i_j\in U\cup V}b^i_j)^2\right)\psi(*)D\Psi_p(\dots)\\
-\sum_U(b^{q_1}_{k_1}+b^{q_2}_{k_2}+1+
\sum_{x^i_j\in U}b^i_j)(\sum_{x^i_j\in
U}b^i_j)D_{U\cup\{x^{q_1}_{k_1}
x^{q_2}_{k_2}\}}D\Psi_p(\dots).
\end{multline}

This expression is obviously symmetric. Hence the right-hand side of
formula~(\ref{double-expression-for-psi-classes}) is symmetric.
\end{proof}

\subsection{Formula for Hurwitz numbers of polynomials}\label{formula-hur-pol}

\begin{theorem}\label{thm1} If $m\geq 3$, then
\begin{equation}\label{formula_1}
h(0,n|A_1,\dots,A_m)
=\frac{n^{m-3}}{\prod_{j=2}^m|\aut(A_j)|} \cdot
\int_{\oM_{0,N}}
\psi(x_1^1)^{m-3}
\prod_{p=2}^m \Psi_p
\end{equation}
\end{theorem}

Recall that by $\Psi_p$ we denote $\Psi_p(b^i_j)$, where $b^i_j=a^i_j-1$,
$i=2,\dots,p-1,p+1,\dots,m$, and $j=1,\dots,l_i$.

Below we check this formula in some special cases independently of its proof. It is probably helpful to look through this section in order to see how to work with classes like $\pi^*_{p,q,i}\psi(x^q_i)$, $D_U$, and $\Psi_p$.




\section{Admissible covers}\label{admissible-covers}

\subsection{Covering of $\M_{0,m}$}

Let us fix the degree $n$, the genus $g$, and a collection of passports $A_1,\dots,A_m$. Here, for any $i=1,\dots,m,$ $A_i=(a^i_1,\dots,a^i_{l_i})$, $a^i_1\geq\dots\geq a^i_{l_i}$, and $\sum_{j=1}^{l_i}a^i_j=n$.

Consider all meromorphic functions $f\colon C\to\CP1$ of degree $n$ defined on smooth curves of genus $g$ such that, for certain distinct points $z_1,\dots,z_m\in\CP1$, the passports of $f$ over $z_1,\dots,z_m$ are equal to $A_1,\dots,A_m$, respectively, and $f$ is unramified over $\CP1\setminus\{z_1,\dots,z_m\}$. We consider such functions up to isomorphism and up to automorphisms of $\CP1$ in the target. Then the space of such function is a noncompact complex manifold of dimension $m-3$. Denote this space by $H$.

Consider $(\CP1,z_1,\dots,z_m)$ as a moduli point of $\M_{0,m}$. Then there is a natural projection $ll\colon H\to\M_{0,m}$. The mapping $ll$ is usually called the \emph{Lyashko-Looijenga mapping}. Obviously, $ll\colon H\to\M_{0,m}$ is an $h(g,n|A_1,\dots,A_m)$-sheeted unramified covering.

\subsection{Boundary of $\oM_{0,m}$}

It order to geometrically obtain a moduli point of $\oM_{0,m}\setminus\M_{0,m}$ with $k$ nodes, one can choose $k$ pairwise nonintersecting contours $c_1,\dots,c_k$ on $(\CP1,z_1,\dots,z_m)$ and contract each of these contours. Contours $c_1, \dots, c_k$ must contain no points $z_1,\dots,z_m$ and no points of self-intersections. The Euler characteristic of each connected component of $\CP1\setminus\left(\bigcup_{i=1}^kc_i\cup\bigcup_{i=1}^mz_i\right)$ must be negative.

We give a more rigorous definition. A moduli point of $\oM_{0,m}$ is a tree of rational curves. Any two irreducible components either are disjoint or intersect transversely at a single point. Each component must contain at least three special (singular or labeled) points. For details see, for example,~\cite{ma,k}.

\subsection{Admissible covers}

Our goal now is to extend the unramified covering $ll\colon H\to\M_{0,m}$ to a ramified covering $\oll\colon\oH\to\oM_{0,m}$.

Suppose that a moduli point on the boundary of $\oM_{0,m}$ is obtained from a moduli point $(\CP1,z_1,\dots,z_m)$ by contracting contours $c_1,\dots,c_k$. Consider a function $f\in H$, $f\colon C\to\CP1$ such that $ll(f)=(\CP1,z_1,\dots,z_m)$. All preimages of contours $c_1,\dots,c_k$ are contours on $C$. We contract each of them.

Thus we obtain a point on the boundary of $\oH$. The axiomatic description of functions we obtain by such procedure gives us the definition of the space $\oH$.

Consider a moduli point $(C,z_1,\dots,z_m)$ on the boundary of $\oM_{0,m}$. Then $\oll^{-1}(C,z_1,\dots,z_m)$ consists of holomorphic maps $f\colon C_g\to C$ of prestable curves with $m$ labeled points to $C$ such that over each irreducible component of $C$ $f$ is an
$n$-sheeted covering, not necessarily connected, with ramifications only over special points (labeled or singular). It is required that ramifications over marked points are determined by their passports, and the local behavior of $f$ at a node in the preimage is the same
on both branches of $C_g$ at this node.

Thus we defined the space $\oH$ and the mapping $\oll$. For more detailed definitions, we refer to~\cite{hm,hmo,i}.

\subsection{Space of admissible covers}

The space $\oH$ can be considered as an orbifold with some glued strata. Actually, it can be desingularized, but we only need this space to carry a fundamental class, see~\cite{av, i}.




\section{Lemma of E.~Ionel}\label{lemma-of-ionel}

\subsection{Covering of space of admissible covers}\label{cov-amiss-covers}

Let us fix the degree $n$, the genus $g$, and a collection of passports $A_1,\dots,A_m$. Consider the corresponding space of admissible covers $\oH$.

We denote by $\hH$ the space of functions from $\oH$ with all labeled preimages of all points $z_1,\dots,z_m$. The natural projection $\pi\colon\hH\to\oH$ is a $\left(\prod_{i=1}^m|\aut(A_i)|\right)$-sheeted ramified covering.

We denote by $\hll\colon\hH\to\oM_{0,m}$ the mapping $\oll\circ\pi$. It is also a ramified covering.

\subsection{Example}

Consider $n=3$, $g=0$, $A_1=(3)$, $A_2=A_3=(2,1)$, and $A_4=(1,1,1)$. Then we expect $\hll\colon\hH\to\oM_{0,4}$ to be a $h(0,3|A_1,\dots,A_4)\cdot \left(\prod_{i=1}^4|\aut(A_i)|\right)=6$-sheeted ramified covering.

Obviously, $\hll$ is ramified only over the boundary points of $\oM_{0,4}$. Denote these points by $(14|23)$, $(12|34)$, and $(13|24)$ (we mean here, that $(14|23)$ corresponds to the curve, where $z_1,z_4$ and $z_2,z_3$ are labeled points on different irreducible components).

It is easy to see that the passport of $\hll$ over $(14|23)$ is equal to $(3,3)$ and the passports of $\hll$ over all other boundary points are equal to $(2,2,2)$. Thus $\hll$ is a $6$-sheeted covering of the sphere and it has $10$ critical points. Hence, the covering space $\hH$ is a sphere.

\subsection{Lemma of E.~Ionel}

Let the passport $A_i$ be equal to $(a^i_1,\dots,a^i_{l_i})$. Then a point of the space $\hH$ is a function
\begin{equation}
f\colon (C_g,x^1_1,\dots,x^1_{l_1},\dots,x^m_1,\dots,x^m_{l_m})\to
(C_0,z_1,\dots,z_m).
\end{equation}
Here $x^i_1,\dots,x^i_{l_i}$ are labeled preimages of $z_i$,
$f^{-1}(z_i)=\sum_{j=1}^{l_i}a^i_jx^i_j$.

By $st\colon\hH\to\oM_{g,N}$, $N=\sum_{i=1}^ml_i$ we denote the mapping that takes a function $\left( f\colon (C_g,x^1_1,\dots,x^m_{l_m})\to (C_0,z_1,\dots,z_m) \right) \in\hH$ to the moduli point $(C_g,x^1_1,\dots,x^1_{l_1},\dots,x^m_1,\dots,x^m_{l_m}) \in\oM_{g,N}$.

\begin{lemma}\cite{i} For any $i=1,\dots,m$ and $j=1,\dots,l_i$
\begin{equation}
a^i_j\cdot st^*\psi(x^i_j)=(\hll)^*\psi(z_i).
\end{equation}
\end{lemma}

Here, by $\psi(x^i_j)$ (by $\psi(z_i)$) we denote the first Chern class of the line bundle whose fiber is the cotangent line at the point $x^i_j$ (resp., at the point $z_i$).




\section{Proof of Theorem~\ref{thm1}}

This section is organized as follows. In Section~\ref{reminder}, we recall notation and the statement of Theorem~\ref{thm1}. In Section~\ref{reminder}, we prove the theorem. This proof is based on the calculation of $st_*[\hH]$ given in Section~\ref{poincare-dual}. In Section~\ref{independent-check}, we check Theorem~\ref{thm1} independently of its proof.

\subsection{Reminder of notation}\label{reminder}

Recall that we fix $n\in\N$, $g=0$, and passports $A_1,\dots,A_m$, $A_i=(a^i_1,\dots,a^i_{l_i})$. We also assume that $l_1=1$, $A_1=(n)$, and $\sum_{i=1}^m\sum_{j=1}^{l_i}(a^i_j-1)=2n-2$.

Then, if $m\geq 3$, we want to prove that
\begin{equation}
h(0,n|A_1,\dots,A_m)
=\frac{n^{m-3}}
{\prod_{j=2}^m |\aut(A_j)|}
\int_{\oM_{0,N}}
\psi(x_1^1)^{m-3}
\prod_{p=2}^m
\Psi_p.
\end{equation}

\subsection{Proof}

Consider the corresponding space $\hH$ defined in Section~\ref{cov-amiss-covers}. We have the following picture:
\begin{equation}(C_0,x^1_1,\dots,x^m_{l_m})\in
\oM_{0,N}\stackrel{st}{\longleftarrow}\hH
\stackrel{\hll}{\longrightarrow}\oM_{0,m}
\ni(C_0,z_1,\dots,z_m)
\end{equation}

Note that $\int_{\oM_{0,m}}\psi(z_1)^{m-3}=1$. Since $\hll$ is a ramified
covering of degree
\begin{equation}
S=h(0,n|A_1,\dots,A_m)\cdot\left(\prod_{i=2}^m|\aut(A_i)|\right),
\end{equation}
it follows that $\int_{\hH}\hll^*\psi(z_1)^{m-3}=S$.

Then, using the lemma of E.~Ionel, we obtain that
\begin{equation}
S=n^{m-3}\int_{\hH}st^*\psi(x^1_1)^{m-3}=
n^{m-3}\int_{st_*[\hH]}\psi(x^1_1)^{m-3}.
\end{equation}

Below we prove that in our case the Poincar\'e dual of $st_*[\hH]$ is equal to $\Xi=\prod_{p=2}^m\Psi_p$.

Thus we obtain our formula.

\subsection{The Poincar\'e dual of $st_*[\hH]$}\label{poincare-dual}

In this subsection we prove that $st_*[\hH]$ is dual to $\Xi$. First we define a subvariety $V\subset\oM_{0,M}$; then we prove that $[V]$ is dual to $\Xi$; and then we prove that $[V]=st_*[\hH]$.

\subsubsection{Subvariety $V$}

Consider a moduli point $(C,x^1_1,\dots,x^m_{l_m})\in\oM_{0,M}$. We fix $p\geq 2$.
Consider the meromorphic $1$-form $\omega_p$ with simple poles at the points $x^1_1,x^p_1,\dots,x^p_{l_p}$ with residues $-n,a^p_1,\dots,a^p_{l_p}$ respectively.

Let $C_p\subset C$ be the union of those irreducible components of $C$, where $\omega_p$ is not identically zero. Obviously, $C_p$ is a connected curve. Consider the collapsing map $c_p\colon C\to C_p$. Since $C_p$ is a connected curve, it follows that the image $c_p(x^i_j)$ of
any labeled point $x^i_j$ is a nonsingular and a nonlabeled point of $(C_p,x^1_1,x^p_1,\dots,x^p_{l_p})$.

Let $x\in C_p$ be the image of the marked points $x^{q_1}_{j_1},\dots,x^{q_{k(x)}}_{j_{k(x)}}$ under the mapping $c_p$. Then we require that $\omega_p$ has exactly $\sum_{i=1}^{k(x)}(a^{q_i}_{j_i}-1)$ zeros at $x$.

We define the subvariety $V$ as follows. A moduli point $(C,x^1_1,\dots,x^m_{l_m})\in\oM_{0,M}$ belongs to $V$ if and only if it satisfies the last requirement for any $p=2,\dots,m$ and for any nonsingular nonlabeled point $x\in (C_p,x^1_1,x^p_1,\dots,x^p_{l_p})$.

\subsubsection{$[V]$ is dual to $\Xi$}

We fix $p\geq 2$ and a sequence of numbers
$(b^i_j)^{i\in\{2,\dots,m\}\setminus\{p\}}_{j=1,\dots,l_i}$. By $V_p(b^i_j)$
denote the subvariety of $\oM_{0,M}$ such that for any
$(C,x^1_1,\dots,x^m_{l_m})\in V$ for any nonsingular nonlabeled point $x\in
C_p$, $x=c_p(x^{q_1}_{j_1})=\dots=c_p(x^{q_{k(x)}}_{j_{k(x)}})$, $\omega_p$ has
at least $\sum_{i=1}^{k(x)}b^{q_i}_{j_i}$ zeros at $x$.

\begin{lemma}\label{v-b-is-dual-to-psi-b} The subvariety $V_p(b^i_j)$ determines
a cohomological class equal to $\Psi_p(b^i_j)$.
\end{lemma}

\begin{proof} Actually, if $b^i_j=0$ for all $i$ and $j$, then
$V_p(b^i_j)=\oM_{0,M}$ and $[V_p(b^i_j)]$ is dual to $\Psi_p(0)=1$.

Suppose that we have already proved this lemma for all classes with
\begin{equation}
\sum_{i\in\{2,\dots,m\}\setminus\{p\}}\sum_{j=1}^{l_i}b^i_j\leq s.
\end{equation}
Consider a sequence $(b^i_j)^{i\in\{2,\dots,m\}\setminus\{p\}}_{j=1,\dots,l_i}$ such that
\begin{equation}
\sum_{i\in\{2,\dots,m\}\setminus\{p\}}\sum_{j=1}^{l_i}b^i_j=s.
\end{equation}
We fix $q\in\{2,\dots,m\}\setminus\{p\}$ and $k\in\{1,\dots,l_q\}$. Then we define $(\widehat b^i_j)^{i=2,\dots,p-1,p+1,\dots,m}_{j=1,\dots,l_i}$ in the following way. We put $\widehat b^q_k=\widehat b^q_k+1$; for all other indices $i,j$ we put $\widehat b^i_j=b^i_j$.

By $\L(x^q_k)$ denote the cotangent line bundle at the point $x^q_k$.

On $V_p(b^i_j)$, the restriction $\omega_p|_{x^q_k}$ determines a section of $\pi^*_{p,q,k}\L(x^q_k)^{\otimes \widehat b^q_k}$. Note that, at a generic point of $V_p(b^i_j)$, the form $\omega_p$ has exactly $b^q_k$ zeros at the point $c_p(x^q_k)$. Hence the restriction $\omega_p|_{x^q_k}$ vanishes at all the curves where $\omega_p$ has more than $b^q_k$ zeros at the point $c_p(x^q_k)$.

Let us enumerate those divisors in $V_p(b^i_j)$, where $\omega_p$ has more than $b^q_k$ zeros at $c_p(x^q_k)$. The first one is the divisor where $x^q_k\in C_p$ and $\omega_p$ has $\widehat b^q_k$ zeros at $x^q_k=c_p(x^q_k)$. The second case is when $x^q_k$ `run away' from $C_p$ with some other points $x^{i_1}_{j_1},\dots,x^{i_s}_{j_s}$. The section of $\pi^*_{p,q,k}\L(x^q_k)^{\otimes \widehat b^q_k}$ has $\sum_{r=1}^sb^{i_r}_{j_r}$ zeros at such divisor.

Hence, the expression defining $\Psi_p(\widehat b^i_j)$ is just the expression defining the same cocycle as the first divisor described above. But the first divisor is just $V_p(\widehat b^i_j)$. Thus we made the inductive step which proves this lemma.
\end{proof}

Note that $V$ is a transversal intersection of $V_p(a^i_j-1)$, $p=2,\dots,m$. Then an obvious corollary of Lemma~\ref{v-b-is-dual-to-psi-b} is the following lemma

\begin{lemma} The subvariety $V$ determines a cohomological class equal to $\Xi$.
\end{lemma}

\subsubsection{$[V]$ is equal to $st_*[\hH]$}

Obviously, $V$ and $st(\hH)$ are the closures of $V\cap\M_{0,M}$ and $st(\hH)\cap\M_{0,M}$,
respectively. We can prove that $V\cap\M_{0,M}=st(\hH)\cap\M_{0,M}$.

Let $f\colon(\CP1,x^1_1,\dots,x^m_{l_m})\to(\CP1,z_1,\dots,z_m)$ be a function in $\hH$. Then
$\omega_p=f^*\left(dz/(z-z_p)-dz/(z-z_1)\right)$. Since $x^q_k$ is a critical point of $f$ of the order $a^q_k-1$, if follows that $\omega_p$ has $a^q_k-1$ zeros at $x^q_k$. Thus we obtain that $V\cap\M_{0,M}\supset st(\hH)\cap\M_{0,M}$.

Conversely, let $(\CP1,x^1_1,\dots,x^m_{l_m})\in V\cap\M_{0,M}$. Then there is a unique function $f_p$ such that $f_p^{-1}(z_1)=n\cdot x^1_1$ and $f_p^{-1}(z_p)=\sum_{i=1}^{l_p} a^p_i\cdot x^p_i$. Requirements for $\omega_p$ mean that a marked point $x^q_k$ is a critical point of $f_p$ of order $a^q_k-1$. Hence $f=f_2=f_3=\dots=f_m$ is the function in $\hH$ such that $f^{-1}(z_i)=\sum_{j=1}^{l_i}a^i_jx^i_j$, $i=1,\dots,m$.

Since the map $st$ has degree one over $V\cap\M_{0,M}$, it follows that $[V]=st_*[\hH]$.




\subsection{Independent check of Theorem~\ref{thm1}}\label{independent-check}

\subsubsection{Degenerate polynomials}

Consider the first nontrivial degenerate case. Let $A_1=(n)$, $A_2=(a^2_1,a^2_2)$, and $A_3=(2,1,\dots,1)$. Then we want to check that $h(0,n|A_1,A_2,A_3)=1/\aut(A_2)$.

We have
\begin{equation}
\begin{array}{rcl}
\Psi_2&=&\pi_{2,3,1}^*\psi(x^3_1),\\
\Psi_3&=&
(a^2_1-1)!(a^2_2-1)!
\pi_{3,2,1}^*\psi(x^2_1)^{a^2_1-1}
\pi_{3,2,2}^*\psi(x^2_2)^{a^2_2-1}\\
&&
-((n-2)!-(a^2_1-1)!(a^2_2-1)!)D_{\{x^2_1,x^2_2\}}
\pi_{3,2,1}^*\psi(x^2_1)^{n-3}.
\end{array}
\end{equation}

Note that $\Psi_2$ is dual to the divisor $D$ whose generic point is represented by a two-component curve such that $x^1_1$ and $x^2_1$ lie on one component and $x^2_2$ and $x^3_1$ lie on the other one.  The restriction of $\Psi_3$ to this divisor is obviously equal to
\begin{equation}
(a^2_1-1)!(a^2_2-1)!
\psi(x^2_1)^{a^2_1-1}\psi(x^2_2)^{a^2_2-1}
\end{equation}

The divisor $D$ consists of several irreducible components. One can enumerate these components via subsets of $\{x^3_2,\dots,x^3_{n-1}\}$. Let $D_U$, $U\subset \{x^3_2,\dots,x^3_{n-1}\}$ be the divisor whose generic point is represented by a two-component curve such that $x^1_1$, $x^2_1$, and all $x^3_i\in U$ lie on one component and all other labeled points lie on the other component. So, $D=\cup_U D_U$.

It is easy to see that $\int_{D_U}\psi(x^2_1)^{a^2_1-1}\psi(x^2_2)^{a^2_2-1}$ is not vanishing if and only if $|U|=a^2_1-1$. In this case this integral is equal to $1$. Hence, $\int_{D}\psi(x^2_1)^{a^2_1-1}\psi(x^2_2)^{a^2_2-1}=\binom{n-2}{a^2_1-1}$. It follows that $h(0,n|A_1,A_2,A_3)=1/\aut(A_2)$.

\subsubsection{Generic polynomials}

The natural desire here is to independently check that out formula gives the same answer as in~\cite{l} in the case of generic polynomials. The ramification data for generic polynomials is the following: $A_1=(n)$, $m=n$, and $A_2=\dots=A_n=(2,1,\dots,1)$. The answer given in~\cite{l} is $h(0,n|A_1,\dots,A_n)=n^{n-3}$.

We state that $h(0,n|A_1,\dots,A_n)$ counted with our formula is equal to $n^{n-3}$. But the argument we have is too complicated, so we give only the sketch of the proof, which can be deciphered to the full proof by the reader.

The first step looks as follows. One has to show that the `essential part' of $\Psi_p$ in this case is equal to the cocycle determined by the subvariety $W_p$. $W_p$ is a codimension $n-2$ subvariety whose generic point is represented by an $(n-2)$-component curve such that on one component there are points $x^1_1$, $x^p_1$, and $n-2$ singular points. All other components, except for this one, contain exactly one point from the set $\{x^2_1,\dots,x^{p-1}_1,x^{p+1}_1,\dots,x^n_1\}$, exactly one point from the set $\{x^p_2,\dots,x^p_n\}$, and exactly one singular point, attaching this component to the first one.

The words `essential part' mean that the integral
\begin{equation}
\int_{\oM_{0,N}}\psi(x^1_1)^{n-3}\Psi_2\dots\Psi_{p-1}(\Psi_p-[W_p])\Psi_{p+1}
\dots\Psi_n
\end{equation}
is equal to zero.

We remark that there are exactly $(n-2)!$ ways to split the sets $\{x^2_1,\dots,x^{p-1}_1,x^{p+1}_1,\dots,x^n_1\}$ and $\{x^p_2,\dots,x^p_n\}$ into pairs of points lying on the same components. Then, via the standard calculations, one can show that the integral part of formula~(\ref{formula_1}) is equal to $(n-2)!^{n-1}$. Then we obtain that
\begin{equation}
h(0,n|A_1,\dots,A_n)=\frac{n^{n-3}}{\prod_{i=1}^n\aut(A_i)}(n-2)!^{n-1}=n^{n-3}.
\end{equation}




\section{Hurwitz numbers of generalized polynomials}

\subsection{Two-pointed ramification cycles}

In applications of the lemma of E.~Ionel to concrete calculations, one has to work with two-pointed ramification cycles. We define them now.

Consider a moduli space $\oM_{g,k}$. Let $x_1,\dots,x_k$ be marked points of curves in $\oM_{g,k}$. Let $b_1,\dots,b_k$ be integer numbers such that $\sum_{i=1}^kb_i=0$. By $W$ denote the subvariety of $\M_{g,k}$ consisting of curves $(C,x_1,\dots,x_k)\in\M_{g,k}$ such that $\sum_{i=1}^kb_kx_k$ is the divisor of a meromorphic function.

The closure of $W$ in $\oM_{g,k}$ determines a homology class $\Delta(b_1,\dots,b_k)$. This class is called the \emph{two-pointed ramification cycle}.

\subsection{Hurwitz numbers of generalized polynomials}\label{subsection-generalized-polynomials}

By a \emph{generalized polynomial} we simply mean the following ramification data. We fix integers
$n>0$ and $g\geq 0$. Let $A_1=(n)$, $A_2=(a_1,\dots,a_l)$, and $A_3=\dots=A_m=(2,1,\dots,1)$, where $m-l-1=2g$ (the Riemann--Hurwitz formula).

Consider the space $\oM_{g,1+l}\ni(C,x^1_1,x^2_1,\dots,x^2_l)$. There is a two-pointed ramification cycle $\Delta(-n,a_1,\dots,a_l)$, where $-n$ corresponds to $x^1_1$, and $a_i$ corresponds to $x^2_i$, $i=1,\dots,l$.

\begin{theorem}\label{theorem-generalized-polynomials} If $m\geq 3$, then
\begin{equation}\label{generalized-polynomials}
h(g,n|A_1,\dots,A_m)=\frac{n^{m-3}(m-2)!}{\aut(A_2)}
\int_{\Delta(-n,a_1,\dots,a_l)} \psi(x_1^1)^{m-3}.
\end{equation}
\end{theorem}

\subsection{Genus zero case}
We check formula~(\ref{generalized-polynomials}) in the case of genus zero. We have $\Delta(-n,a_1,\dots,a_l)=[\oM_{0,1+l}]$. Hence,
\begin{equation}
h(0,n|A_1,\dots,A_{l+1})=\frac{n^{l-2}(l-1)!}{\aut(A_2)}.
\end{equation}
The same answer is given by~(\ref{formula-gjlz}).




\section{Proof of Theorem~\ref{theorem-generalized-polynomials}}

\subsection{The first steps of the proof}

We start our proof in the same way as in the case of usual polynomials. Let $\oH$ be the appropriated space of admissible covers. The mapping $\hH\to\oH$ has degree $\aut(A_2)\cdot ({n-2})!^{m-2}$. Then,
\begin{equation}
\int_{\hH}\hll^*\psi(z_1)^{m-3}=\aut(A_2)\cdot (n-2)!^{m-2}\cdot
h(g,n|A_1,\dots,A_m).
\end{equation}
Using the lemma of E.~Ionel, we get
\begin{equation}
n^{m-3}\int_{st_*[\hH]}\psi(x^1_1)^{m-3}=
\aut(A_2)\cdot (n-2)!^{m-2}\cdot
h(g,n|A_1,\dots,A_m).
\end{equation}

\subsection{The restriction of $\psi(x^1_1)$ to $st(\hH)$}

Let $\pi\colon\oM_{g,1+l+(m-2)(n-1)}\to\oM_{g,1+l}$ be the projection forgetting
all marked points except for $x^1_1,x^2_1,$ $\dots,x^2_l$.

\begin{lemma}\label{lemma-shifted-psi}
$\pi^*\psi(x^1_1)|_{st(\hH)}=\psi(x^1_1)|_{st(\hH)}$. \end{lemma}

\begin{proof}
Recall that $\pi^*\psi(x^1_1)=\psi(x^1_1)-[D_\pi]$. Here $[D_\pi]$ is the
cocycle determined by the divisor $D_\pi$ in $\oM_{g,1+l+(m-2)(n-1)}$. The
generic point of $D_\pi$ is represented by a two-component curve such that one
component has genus zero and contains $x^1_1$, and the other component has genus
$g$ and contains $x^2_1,\dots,x^2_l$.

Let us prove that this divisor does not intersect $st(\hH)$. Assume the converse, that is, we assume that there exists a function in $\hH$ such that its domain belongs to $D_\pi$.

Consider the image of such function. It is a stable curve of genus zero with marked points $z_1,\dots,z_m$. It follows from the definition of $D_\pi$ that $z_1$ and $z_2$ lie on different irreducible components of the target curve.

Consider the irreducible component $C_1$ of the target curve containing $z_1$. This component also contains the point $z_*$ and at least one more special point. (The point $z_*$ separates $C_1$ from the component containing $z_2$.)

If follows from the definition of admissible covers that the point $z_*$ is a point of total ramification. Since the preimage of $C_1$ has genus zero, it follows that any point of this component (except for $z_1$ and $z_*$) has exactly $n$ simple preimages.

Note that all marked points in the target curve are critical values of the function. After being cut at $z_*$, the target curve splits in two halves, and the total preimage of the half containing $z_1$ has genus $0$. It follows that the third special point on $C_1$ has a ramification in the preimage.

This contradiction proves that $D_\pi$ does not intersect $st(\hH)$. Hence, $\pi^*\psi(x^1_1)|_{st(\hH)}=\psi(x^1_1)|_{st(\hH)}$.
\end{proof}

\subsection{Proof of Theorem~\ref{theorem-generalized-polynomials}}

We have
\begin{equation}
n^{m-3}\int_{st_*[\hH]}\psi(x^1_1)^{m-3}=
\aut(A_2)\cdot (n-2)!^{m-2}\cdot
h(g,n|A_1,\dots,A_m).
\end{equation}
Using Lemma~\ref{lemma-shifted-psi}, we get
\begin{equation}
\int_{st_*[\hH]}\psi(x^1_1)^{m-3}=\int_{\pi_*st_*[\hH]}\psi(x^1_1)^{m-3}.
\end{equation}
Since the mapping $\pi\circ st|_\hH$ has degree $(m-2)!\cdot (n-2)!^{(m-2)}$, it follows that
\begin{equation}
\pi_*st_*[\hH]=(m-2)!\cdot (n-2)!^{(m-2)}\cdot\Delta(-n,a_1,\dots,a_l).
\end{equation}
Hence,
\begin{equation}
h(g,n|A_1,\dots,A_m) =
\frac{n^{m-3}\cdot (m-2)!}{\aut(A_2)}\cdot
\int_{\Delta(-n,a_1,\dots,a_l)}\psi(x^1_1)^{m-3}.
\end{equation}




\section{Amusing formulas for $\l\tau_{3g}\tau_0^2\r_g$}

Consider the moduli space of curves $\oM_{g,n}\ni(C,x_1,\dots,x_n)$. By $\l\tau_{k_1}\dots\tau_{k_n}\r_g$ denote $\int_{\oM_{g,n}}\psi(x_1)^{k_1}\dots\psi(x_n)^{k_n}$. It is known from~\cite{w,ko} that $\l\tau_{3g}\tau_0^2\r_g=1/(24^gg!)$.

In this section, we express $\l\tau_{3g}\tau_0^2\r_g$ via Hurwitz numbers of generalized polynomials.

\subsection{Formulas}\label{formulas-tau3g}

By $H(g;n)$ we denote $h(g,n|A_1,\dots,A_m)$, where
$A_1=(n)$, $A_2=\dots=A_m=(2,1,\dots,1)$, $m=2g+n$.

\begin{theorem}\label{formulas-for-tau3g}
For any $l\geq 0$,
\begin{equation}\label{formula-tau3g}
\l\tau_{3g}\tau_0^2\r_g=\sum_{i=0}^{g}
(-1)^i\binom{g}{i}
\frac{(g+l+1-i)!H(g;g+l+1-i)}
{g!(3g+l-i)!(g+l+1-i)^{3g+l-1-i}}.
\end{equation}
\end{theorem}

Thus we have an infinite number (for any $l\geq 0$) of formulas for $\l\tau_{3g}\tau_0^2\r_g$.

\subsection{Check in low genera}

\subsubsection{Genus zero} In genus zero, we have the following:
\begin{equation}
\l\tau_0^3\r_0=\frac{(l+1)!H(0;l+1)}{l!(l+1)^{l-1}}.
\end{equation}
Since $H(0;l+1)=(l+1)^{l-2}$, if follows that
\begin{equation}
\l\tau_0^3\r_0=\frac{(l+1)!(l+1)^{l-2}}{l!(l+1)^{l-1}}=1.
\end{equation}

\subsubsection{Genus one} In genus one, we have the following:
\begin{equation}
\l\tau_3\tau_0^2\r_1=
\frac{(l+2)!H(1;l+2)}
{(l+3)!(l+2)^{l+2}}
-
\frac{(l+1)!H(1;l+1)}
{(l+2)!(l+1)^{l+1}}.
\end{equation}
From~\cite{v}, it follows that
\begin{equation}
H(1;l)=\frac{(l+1)!l^l(l-1)}{l!\cdot 24}.
\end{equation}
Combining these formulas, we obtain that $\l\tau_3\tau_0^2\r_1=1/24$.




\section{Proofs of Theorem~\ref{formulas-for-tau3g}}

We will give two proofs of Theorem~\ref{formulas-for-tau3g}. The first one is purely combinatorial. It is based on the formula of Ekedahl et al.~\cite{elsv}. The second proof is purely geometric. This proof is based on Theorem~\ref{theorem-generalized-polynomials}.

Historically, the second proof was the first one. It is based on ideas which will be very useful in the rest of the paper. We obtained the second proof trying to check our formula in genera $2$ and $3$.

\subsection{First proof}

\begin{proof} We recall the formula for $H(g;k)$ from~\cite{elsv}
\begin{equation}
H(g;k)=\frac{(2g+k-1)!\cdot
k^k}{k!}\int_{\oM_{g,1}}
\frac{\sum_{i=0}^g(-1)^i\lambda_i}
{1-k\psi(x_1)}.
\end{equation}
Here $x_1$ is the unique marked point on curves in $\oM_{g,1}$, and $\lambda_i=c_i(\E)$, where $\E$ is the rank $g$ vector bundle over $\oM_{g,1}$ with the fiber over $(C,x_1)$ equal to $H^0(C,\omega_C)$.

Putting this formula for $H(g;k)$ in formula~(\ref{formula-tau3g}), we obtain the following:
\begin{equation}
\l\tau_{3g}\tau_0^2\r_g=\frac{1}{g!}
\sum_{i=0}^{g}
\sum_{j=0}^{g}
(-1)^{i+j}\binom{g}{j}(g+l+1-j)^{g-i}
\int_{\oM_{g,1}}\lambda_i\psi(x_1)^{3g-2-i}.
\end{equation}

It follows from the combinatorial lemma below that the right hand side of this formula is equal to $\int_{\oM_{g,1}}\psi(x_1)^{3g-2}$. The equality $\l\tau_{3g}\tau_0^2\r_g=\l\tau_{3g-2}\r_g$ follows from the string equation~\cite{w}:
\begin{equation}
\l\tau_0\prod_{i=1}^m\tau_{k_i}\r=\sum_{i=1^m}\l\tau_{k_1}\dots\tau_{k_{i-1}}
\tau_{k_i-1}\tau_{k_{i+1}}\dots\tau_{k_m}\r.
\end{equation}
\end{proof}

\begin{lemma} There is a combinatorial identity:
\begin{equation}
\sum_{i=0}^g(-1)^i\binom{g}{i}(g+1-i)^k=
\left\{\begin{array}{ll}
0, & k<g, \\
g!, & k=g.
\end{array}\right.
\end{equation}
\end{lemma}

\begin{proof}
By $f(g,k)$ denote the left-hand side of the lemma statement. Note that $f(0,0)=1$, and $f(g,0)=0$ if $g>0$. Since
\begin{equation}
f(g+1,k+1)=
(g+1)\cdot\left(
\binom{k+1}{0}f(g,k)+\dots+\binom{k+1}{k}f(g,0)
\right),
\end{equation}
the statement of the lemma follows.
\end{proof}

\subsection{Second proof}

First, we express the intersection number $\l\tau_{3g}\tau_0^2\r_g=\l\tau_{3g-2}\r_g$ via integrals over two-pointed ramification cycles. Then, using Theorem~\ref{theorem-generalized-polynomials}, we obtain formula~(\ref{formula-tau3g}).

\subsubsection{Two-pointed ramification cycles}

We fix $g\geq 0$. By $\Delta_k$ denote the two-pointed ramification cycle $\Delta(-k,1,\dots,1)$ in the moduli space $\oM_{g,1+k}\ni(C,y,t_1,\dots,t_k)$.

Consider the intersection number $\l\tau_{3g}\tau_0^{2}\r_g$. We prove the following lemma.

\begin{lemma}\label{lemma-proof2}
For any $l\geq 0$,
\begin{equation}
\l\tau_{3g}\tau_0^{2}\r_g=\sum_{i=0}^g
\frac{(-1)^g}{g!}\binom{g}{i}
\int_{\Delta_{g+l-i+1}}\psi(y_1)^{3g+l-i-1}.
\end{equation}
\end{lemma}

\subsubsection{Proof of Theorem~\ref{formulas-for-tau3g}}

\begin{proof}
From Theorem~\ref{theorem-generalized-polynomials}, we know that
\begin{equation}
\int_{\Delta_{g+l-i+1}}\psi(y_1)^{3g+l-i-1}=
\frac{(g+l-i+1)!H(g,g+l-i+1)}{(g+l-i+1)^{3g+l-i-1}(3g+l-i)!}.
\end{equation}
If we combine this with Lemma~\ref{lemma-proof2}, we get
\begin{equation}
\l\tau_{3g}\tau_0^{2}\r_g=\sum_{i=0}^g
\frac{(-1)^g}{g!}\binom{g}{i}
\frac{(g+l-i+1)!H(g,g+l-i+1)}{(g+l-i+1)^{3g+l-i-1}(3g+l-i)!}.
\end{equation}
\end{proof}

\subsubsection{Proof of Lemma~\ref{lemma-proof2}}

Consider the moduli space $\oM_{g,2+g+l}\ni(C,y,$ $t_1,\dots,t_{g+l+1})$. By $V(k_1|i_1,\dots,i_{k_2})$ denote the subvariety of $\M_{g,2+g+l}$ consisting of curves $(C,y_1,t_1,\dots,t_{g+l+1})$ such that there exists a meromorphic function of degree $k$ with pole of multiplicity $k$ at $y$ and simple zeros at $t_{i_1},\dots,t_{i_{k_2}}$.

Let $E(k_1|i_1,\dots,i_{k_2})$ be the cycle in homologies of $\oM_{g,1+g+l+1}$
determined by the closure of $V(k_1|i_1,\dots,i_{k_2})$. By $S(k_1,k_2)$ denote
\begin{equation}
\int_{E(k_1|i_1,\dots,i_{k_2})}\psi(y)^{3g+l-1+k_1-k_2}.
\end{equation}
Obviously, the last number depends only on $k_1$ and $k_2$, but not on the choice of $i_1,\dots,i_{k_2}$.

Below, we always suppose that $k_1\leq k_2+g$ and $k_2\geq 1$.

The proof of Lemma~\ref{lemma-proof2} is based on the following lemma.

\begin{lemma}\label{S-expression-via-difference}
If $k_1>k_2$, then
\begin{equation}
S(k_1,k_2)=\frac{1}{k_2-k_1}\left(S(k_1,k_2+1)-S(k_1-1,k_2)\right).
\end{equation}
\end{lemma}

The dependence of $S(k_1,k_2)$ on $l$ is explained by the following lemma.

\begin{lemma}\label{lemma-l-independence}
The number $S(k_1,k_2)$ does not depend on the choice of $l$ that satisfies $g+l+1\geq k_2$.
\end{lemma}

We prove Lemmas~\ref{S-expression-via-difference} and~\ref{lemma-l-independence} in the next subsection.

\begin{proof}[Proof of Lemma~\ref{lemma-proof2}]

Note that $[\oM_{0,2+g+l}]=E(g+l+1|1,2,\dots,l+1)$. (On the generic curve
$(C,y,t_1,\dots,t_{l+1})$, there exists a meromorphic function of degree $g+l+1$
such that $y$ is the pole of multiplicity $g+l+1$ and $t_1,\dots,t_{l+1}$ are
simple zeros. For the proof of such statements, see~\cite[Section 7]{m}). Then,
\begin{equation}
\l\tau_{3g}\tau_0^{2}\r_g=\int_{\oM_{g,g+l+2}}\psi(y)^{4g+l-1}=S(g+l+1|l+1).
\end{equation}

Lemma~\ref{S-expression-via-difference} implies that
\begin{multline}\label{formula-s-1}
S(g+l+1|l+1)=\frac{1}{g}S(g+l+1|l+2)-\frac{1}{g}S(g+l|l+1)\\
=\frac{S(g+l+1|l+3)}{g(g-1)}-\frac{2\cdot S(g+l|l+2)}{g(g-1)}+
\frac{S(g+l-1|l+1)}{g(g-1)}\\
=\dots=
\sum_{i=0}^g\frac{(-1)^i}{g!}\binom{g}{i}S(g+l-i+1|g+l-i+1).
\end{multline}

Lemma~\ref{lemma-l-independence} implies that
\begin{equation}\label{formula-s-2}
S(g+l-i+1|g+l-i+1)=\int_{\Delta_{g+l-i+1}}\psi(y)^{3g+l-i-1}.
\end{equation}

Combining (\ref{formula-s-1}) and (\ref{formula-s-2}) we obtain the statement
of Lemma~\ref{lemma-proof2}.
\end{proof}

\subsubsection{Proofs of Lemmas~\ref{S-expression-via-difference}
and~\ref{lemma-l-independence}}

\begin{proof}[Proof of Lemma~\ref{lemma-l-independence}]
Let $E=E(k_1|1,\dots,k_2)$ be the cycle in the homology of $\oM_{g,2+g+l}$. Let $\widehat E=E(k_1|1,\dots,k_2)$ be the cycle in the homology of $\oM_{g,2+g+\widehat l}$, where $\widehat l=l+1$.

By $\pi\colon\oM_{g,2+g+\widehat l}\to\oM_{g,2+g+l}$ denote the projection forgetting the labeled point $t_{g+\widehat l+1}$. Note that $\psi(y)$ in the cohomology of $\oM_{g,2+g+\widehat l}$ is equal to $\pi^*\psi(y)+D$, where $D$ is the class dual to the divisor whose generic point is represented by a two-component curve such that one component has genus zero and contains points $y$ and $t_{g+\widehat l+1}$, and the other component has genus $g$ and contains all labeled points except for $y$ and $t_{g+\widehat l+1}$.

It is easy to see that $(\pi^*\psi(y)+D)^{K}=\pi^*\psi(y)^K+D\cdot\pi^*\psi(y)^{K-1}$. Also note that $\pi_*(\widehat E\cdot D)=E$.

Since (from dimensional conditions)
\begin{equation}
\int_{\widehat E}\pi^*\psi(y)^{3g+\widehat l-1+k_1-k_2}=0,
\end{equation}
it follows that
\begin{equation}
\int_{\widehat E}\psi(y)^{3g+\widehat l-1+k_1-k_2}=
\int_{\widehat E}D\cdot \pi^*\psi(y)^{3g+l-1+k_1-k_2}=
\int_{E}\psi(y)^{3g+l-1+k_1-k_2}.
\end{equation}

This proves the Lemma.
\end{proof}

\begin{proof}[Proof of Lemma~\ref{S-expression-via-difference}]
Using Lemma~\ref{lemma-l-independence} we can consider $S(k_1,k_2)$ as
\begin{equation}
\int_{E(k_1|1,\dots,k_2)}\psi(y)^{2g+k_1-1},
\end{equation}
where $g+l=k_2$ ($l$ can be negative).

Consider the following ramification data: $A_1=(k_1)$, $A_2=\dots=A_{2g+k_1}=(2,1,\dots,1)$, $A_{2g+k_1+1}=A_{2g+k_1+2}=(1,\dots,1)$. Consider the space $\hH$ built using this ramification data. The map $st$ takes a point of $\hH$ to a curve
\begin{multline}
(C,x^1_1,x^2_1,\dots,x^2_{k_1-1},\dots,x^{2g+k_1}_1,\dots,x^{2g+k_1}_{k_1-1},
x^{2g+k_1+1}_1,\dots,\\
x^{2g+k_1+1}_{k_1},
x^{2g+k_1+2}_1,\dots,x^{2g+k_1+2}_{k_1})\in\oM_{g,1+2k_1+(2g+k_1-1)(k_1-1)}.
\end{multline}

Consider the projection $\pi\colon\oM_{g,1+2k_1+(2g+k_1-1)(k_1-1)}\to\oM_{g,2+k_2}$, which takes the labeled point $x^1_1$ to $y$; $x^{2g+k_1+1}_1,\dots,x^{2g+k_1+1}_{k_2}$ to $t_1,\dots,t_{k_2}$; $x^{2g+k_1+2}_1$ to $t_{k_2+1}$; and forgets all other labeled points.
Note that $\pi_*st_*[\hH]=K\cdot E(k_1|1,\dots,k_2)$, where
\begin{equation}
K=(k_1-1)!(k_1-k_2)!(2g+k_1-1)!(k_1-2)^{(2g+k_1-1)}.
\end{equation}

The Lemma of E.~Ionel implies that
\begin{equation}
\psi(y)\cdot
E(k_1|1,\dots,k_2)=\frac{1}{K}
\pi_*\left(\frac{1}{k_1}st_*(\hll^*\psi(z_1) \cdot [\hH])-D\cdot
st_*[\hH]\right),
\end{equation}
where $D$ is dual to the divisor whose generic point is represented by a two-component curve such that one component has genus $g$ and contains the preimages of the points $t_1,\dots,t_{k_2+1}$, and the other component has genus zero and contains the point $x^1_1$. Here we use the standard expression for $\psi(x^1_1)$ via $\pi^*\psi(y)$.

One can consider $D\cdot st_*[\hH]$ as the class determined by $D\cap st(\hH)$. We are interested only in the irreducible components of $D\cap st(\hH)$ whose image under the map $\pi$ has codimension one.

We describe the generic point of such component. It is a three-component curve. One component has genus zero and contains $x^1_1$, the whole tuple of preimages of one critical value, say $x^2_1,\dots,x^2_{k_1-1}$, and two nodes, $*_1$ and $*_2$. There exists a meromorphic function of degree $k_1$ whose divisor is $-k_1x^1_1+(k_1-1)*_1+*_2$, and $x^2_1,\dots,x^2_{k_1-1}$ are the preimages of its simple critical value. The second component has genus zero and is attached to the first one at the point $*_2$. It contains exactly one point from each of the sets $\{x^3_2,\dots,x^3_{k_1-1}\}$, $\dots$, $\{x^{2g+k_1}_2,\dots,x^{2g+k_1}_{k_1-1}\}$, $\{x^{2g+k_1+1}_{k_2+1},\dots,x^{2g+k_1+1}_{k_1}\}$,
and $\{x^{2g+k_1+2}_{2},\dots,x^{2g+k_1+2}_{k_1}\}$. The last component has genus $g$, is attached to the first component at the point $*_1$, and contains all other points. There is a function of degree $k_1-1$ such that $*_1$ is the point of total ramification, $x^3_1,\dots,x^{2g+k_1}_1$ are simple critical points, and all points $x^i_j$ with fixed superscript are the preimages of a single point in the image.

From this description, it is obvious that
\begin{equation}
\pi_*(D\cdot st_*[\hH])=K\cdot
E(k_1-1|1,\dots,k_2).
\end{equation}

Now we will describe $\pi_*st_*(\hll^*\psi(z_1))$. Note that $\psi(z_1)$ is dual to the divisor whose generic point is represented by a two-component curve such that $z_1$ lies on one component, and $z_{2g+k_1+1}$ and $z_{2g+k_1+2}$ lie on the other component. We describe only those components of the preimage of this divisor under the map $\hll$ whose $\pi_*st_*$-image does not vanish $\psi(y)^{2g+k_1-1}$.

There are two possible cases. The first one looks as follows. Consider the divisor in $\oM_{0,2g+k_1+2}$ whose generic point is represented by a two-component curve such that $z_1$ and exactly one point from the set $\{z_2,\dots,z_{2g+k_1}\}$ lie on one component and all other points (including $z_{2g+k_1+1}$ and $z_{2g+k_1+2}$) lie on the other component. Then we obtain just the same picture as in the case of $D\cdot st_*[\hH]$, but with coefficient $k_1-1$, since $\hll$ is ramified along this divisor with multiplicity $k_1-1$.

The next case is as follows. Consider the divisor in $\oM_{0,2g+k_1+2}$ whose generic point is represented by a two-component curve such that $z_{2g+k_1+1}$ and $z_{2g+k_1+2}$ lie on one component and all other marked points lie on the other component. The map $\hll$ is not ramified over this divisor. There are $k_2+(k_1-k_2)$ possible cases in the preimage of this divisor under the map $\hll$.

The first $k_2$ cases mean that $x^{2g+k_1+2}_1$ lies on the same genus zero irreducible component of the domain curve as one of the marked points $x^{2g+k_1+1}_{1},\dots,x^{2g+k_1+1}_{k_2}$. Under the map $\pi_*st_*$ this gives us $k_2\cdot K\cdot E(k_1|1,\dots,k_2)$ in the space $\oM_{g,2+l'+g}$, where $l'=l-1$.

The other $k_1-k_2$ cases mean that $x^{2g+k_1+2}_1$ lie on the same genus zero irreducible component of the domain curve as one of the marked points $x^{2g+k_1+1}_{k_2+1},\dots,x^{2g+k_1+1}_{k_1}$. Under the map $\pi_*st_*$ this gives us $K\cdot E(k_1|1,\dots,k_2,k_2+1)$ in the space $\oM_{g,2+l+g}$ (it is easy to check the coefficient by direct calculations).

Thus we obtain that
\begin{multline}
\int_{E(k_1|1,\dots,k_2)}\psi(y)^{2g+k_1-1}=
\frac{k_1-1}{k_1}\int_{E(k_1-1|1,\dots,k_2)}\psi(y)^{2g+k_1-2}\\
+\frac{k_2}{k_1}\int^*_{E(k_1|1,\dots,k_2)}\psi(y)^{2g+k_1-2}
+\frac{1}{k_1}\int_{E(k_1|1,\dots,k_2,k_2+1)}\psi(y)^{2g+k_1-2}\\
-\int_{E(k_1-1|1,\dots,k_2,k_2)}\psi(y)^{2g+k_1-2}.
\end{multline}
Here the sign $\int^*$ means that we calculate the intersection number in the (co)homology of $\oM_{g,2+l'+g}$.

Combining the last equality and Lemma~\ref{lemma-l-independence}, we obtain that
$S(k_1,k_2)=(S(k_1,k_2+1)-S(k_1-1,k_2))/(k_1-k_2)$.
\end{proof}




\section{The conjecture of E.~Witten}

In this section, we introduce intersection numbers which will be considered in the rest of the paper. There is a conjecture of E.~Witten relating these intersection numbers with the string solutions of the Gelfand--Dikii hierarchies. Here we also recall this conjecture.

In Section~\ref{inter-num}, we give a definition of the intersection numbers. It is enough to read only this section to understand the rest of the paper except for the appendix.

In Section~\ref{wittens}, we formulate the Witten's conjecture. We recall the necessary definitions for this conjecture in Section~\ref{gel-dikii}. In Section~\ref{relation}, we give a recursive relation for the coefficients of the string solution of the Boussinesq hierachy. This could be considered as an example of the definitions of Section~\ref{gel-dikii}, but we also use this relation in the appendix.

We explain one more time the role of the Witten's conjecture in this paper. There are certain intersection numbers defined by Witten. We give a way to calculate some of these intersection numbers. This is all what we do in the rest of the paper.

But there is a possible application of these results. One can calculate an intersection number using our algorithm. Then one can calculate the corresponding coefficient of the string solution of the corresponding Gelfand--Dikii hierarhy. After this one can compare if the results of these calculations coinside. If yes, then we have checked the Witten conjecture in a very particular case. If no, then the Witten conjecture is false. This is what we do in the appendix for the intersection number $\l\tau_{3,1}\r_3$ in the case of the Boussinesq hierachy.

More detailed exposition of the Witten conjecture can be found in~\cite{w1, jkv}.

\subsection{Intersection numbers
$\l\prod\tau_{i,j}^{k_{i,j}}\r_g$}\label{inter-num}

Some of our definitions in this section seem to be not entirely clear, but our goal now is only to give a general idea of Witten's definitions. As usual, we refer to \cite{w1,jkv,pv} for details.

\subsubsection{A covering of $\M_{g,s}$}

Consider the moduli space $\oM_{g,s}\ni(C,x_1,\dots,x_s)$. Fix an integer $r\geq 2$. Label each marked point $x_i$ by an integer $m_i$, $0\leq m_i\leq r-1$.

By $K$ denote the canonical line bundle of $C$. Consider the line bundle $S=K\otimes(\otimes_{i=1}^s\O(x_i)^{-m_i})$ over $C$. If $2g-2-\sum_{i=1}^sm_i$ is divisible by $r$, then there are $r^{2g}$ isomorphism classes of line bundles $\T$ such that $\T^{\otimes r}\cong S$.

The choice of an isomorphism class of $\T$ determines a cover $\M'_{g,s}$ of $\M_{g,s}$. To extend it to a covering of $\oM_{g,s}$ we have to discuss the behavior of $\T$ near a double point.

\subsubsection{Behavior near a double point}

Let $C$ be a singular curve with one double point. By $\pi\colon C_0\to C$ denote its normalization. The preimage of the double point consists of two points, say $x'$ and $x''$. There are $r$ possible cases of behavior of $\T$ near the double point.

Cases $1,2,\dots,r-1$. The first $r-1$ cases are the following ones: $\T\cong\pi_*\T'$, where $\T'$ is a locally free sheaf on $C_0$ with a natural isomorphism ${\T'}^{\otimes r}\cong
K\otimes(\otimes_{i=1}^s\O(x_i)^{-m_i}) \otimes\O(x')^{-m}\otimes\O(x'')^{-(r-2-m)}$, $m=0,\dots,r-2$.

Case $r$. The last case is when $\T$ is defined by the following exact sequence: $0\to\T'\to\T\to\O\to$. Here $\T'=\pi_*\T''$, where ${\T''}^{\otimes r}\cong
K\otimes(\otimes_{i=1}^s\O(x_i)^{-m_i}) \otimes\O(x')^{-(r-1)}\otimes\O(x'')^{-(r-1)}$. The map $\T\to\O$ is the residue map taking a section of $\T$ to the coefficient of $(dx/x)^{1/r}$.

\subsubsection{The top Chern class}

If $H^0(C,\T)$ vanishes everywhere, then one can consider the vector bundle $\mathcal{V}$ over $\oM'_{g,s}$ whose fiber is the dual space to $H^1(C,\T)$. What we need is the top Chern class
$c_D(\mathcal{V})$ of this bundle; here
\begin{equation}
D=\frac{(g-1)(r-2)}{r}+\frac{1}{r}\sum_{i=1}^sm_i.
\end{equation}

In the case when $H^0(C,\T)$ is not identically zero, there is another definition of the corresponding cohomology class which we will not recall here. We denote this class by the same notation, $c_D(\mathcal{V})$.

We will use only the following properties of this class. First, consider a component of the boundary consisting of curves with one double point, where $\T$ is defined by one of the first $r-1$ possible cases. If this component of the boundary consists of two-component curves, then $c_D(\V)=c_{D_1}(\V_1)\cdot c_{D_2}(\V_2)$, where $\V_1$ and $\V_2$ are the corresponding spaces on the components, $\V=\V_1\oplus\V_2$. If this component of the boundary consists of one-component self-intersecting curves, then $c_D(\V)=c_D(\V')$, where $c_D(\V')$ is defined on the normalization of these curves. Second, consider a component of the boundary consisting of curves with one double point, where $\T$ is defined by the $r$th case. In this case, we require that the rescriction of $c_D(\V)$ to this component of the boundary vanishes.

These properties uniquely determine the desired intersection numbers in the particular case considered below. Nevertheless, it is not obvious that there is a cohomological class which satisfies such properties. As far as we understand, the last problem is clarified in~\cite{jkv,pv}.

\subsubsection{The Mumford-Morita-Miller intersection numbers}

Let us label each marked point $x_i$ by an integer $n_i\geq 0$. By
\begin{equation}
\l\tau_{n_1,m_1}\dots\tau_{n_s,m_s}\r_g
\end{equation}
denote the intersection number
\begin{equation}
\frac{1}{r^g}\int_{\oM'_{g,s}}\prod_{i=1}^s\psi(x_i)^{n_i}\cdot
 c_D(\mathcal{V}).
\end{equation}

Of course, this number is not zero only if
\begin{equation}
3g-3+s=\sum_{i=1}^sn_i+D.
\end{equation}

\subsection{Witten's conjecture}\label{wittens}
Consider the formal series $F$ in variables $t_{n,m}$, $n=0,1,2,\dots$; $m=0,\dots,r-1$;
\begin{equation}
F(t_{0,0},t_{0,1},\dots)=\sum_{d_{n,m}}\l\prod_{n,m}\tau_{n,m}^{d_{n,m}}\r
\prod_{n,m}\frac{t_{n,m}^{d_{n,m}}}{d_{n,m}!}.
\end{equation}

The conjecture is that this $F$ is the string solution of the $r$-Gelfand--Dikii hierarchy.

\subsection{Gelfand--Dikii hierarchies}\label{gel-dikii}
We fix $r$. In this section, we define the string solution of the $r$-Gelfand--Dikii (or $r$-KdV) hierarchy.

Consider the differential operator
\begin{equation}
Q=D^r-\sum_{i=0}^{r-2}u_i(x)D^i,
\end{equation}
where
\begin{equation}
D=\frac{i}{\sqrt{r}}\frac{\d}{\d x}.
\end{equation}

There is a pseudo-differential operator $Q^{1/r}=D+\sum_{i>0}w_iD^{-i}$. By $Q^{n+m/r}$ denote $(Q^{1/r})^{nr+m}$. Note that $[Q^{n+m/r}_+,Q]$ is a differential operator of order at most $r-2$.

The Gelfand--Dikii equations read
\begin{equation}
i\frac{\d Q}{\d t_{n,m}}
=
[Q^{n+m/r}_+,Q]\cdot \frac{c_{n,m}}{\sqrt{r}},
\end{equation}
where
\begin{equation}
c_{n,m}=\frac{(-1)^{n}r^{n+1}}
{(m+1)(r+m+1)\dots(nr+m+1)}.
\end{equation}

The string solution of the Gelfand--Dikii hierarchy is the formal series $F$ in
variables $t_{i,j}$, $i=0,1,2,\dots$, $j=0,\dots,r-1$, such that
\begin{equation}
\frac{\d F}{\d t_{0,0}}=\frac{1}{2}
\sum_{i,j=0}^{r-2}\delta_{i+j,r-2}t_{0,i}t_{0,j}+
\sum_{n=1}^{\infty}
\sum_{m=0}^{r-2}
t_{n+1,m}
\frac{\d F}{\d t_{n,m}},
\end{equation}
\begin{equation}
\frac{\d ^2 F}{\d t_{0,0} \d t_{n,m}}
=-c_{n,m}\res(Q^{n+\frac{m+1}{r}}),
\end{equation}
where $Q$ satisfies the Gelfand--Dikii equations and $t_{0,0}$ is identified with $x$.

One can prove that $F$ is uniquely determined by this equations, up to an additive constant. Below, we will discuss an effective way to calculate the coefficients of $F$ in the case of $r=3$ (Boussinesq hierarchy).

\subsection{Boussinesq hierarchy}\label{relation}

We need to calculate the coefficients of $F$ starting from the Gelfand--Dikii hierarchy. The general methods, like~\cite{n}, are very complicated. We show how to do this in the case $r=3$ (Boussinesq hierarchy).

We denote $\d^k F/\d t_{i_1,j_1}\dots\d t_{i_k,j_k}$ by $\ll \tau_{i_1,j_1}\dots\tau_{i_k,j_k}\rr$. Note that
\begin{equation}
\res(Q^{1/3})=-\frac{\ll\tau_{0,0}^2\rr}{3}, \qquad
\res(Q^{2/3})=-\frac{2\ll\tau_{0,0}\tau_{0,1}\rr}{3}.
\end{equation}

Since $Q=D^3+3(\res(Q^{1/3}))D+(3/2)(\res(Q^{2/3})+D\,\res(Q^{1/3}))$, it follows that $Q=D^3+\gamma_1D+\gamma_2$, where $\gamma_1=-\ll\tau_{0,0}^2\rr$ and $\gamma_2=-\ll\tau_{0,0}\tau_{0,1}\rr-i\ll\tau_{0,0}^3\rr/(2\sqrt{3})$.

We consider the pseudodifferential operator $Q^{n-1+m/r}$. Since $[Q,Q^{n-1+m/r}]=0$, it follows  that $[Q,Q^{n-1+m/r}_-]$ is a differential operator. Let
\begin{equation}
Q^{n-1+m/r}_-=\alpha_1D^{-1}+\alpha_2D^{-2}+\dots.
\end{equation}
Then the coefficients at $D^{-1}$ and $D^{-2}$ of the operator $[Q,Q^{n-1+m/r}_-]$ are
equal to
\begin{equation}\label{bouss-eq-1}
D^3\alpha_1+3D^2\alpha_2+3D\alpha_3+D(\alpha_1\gamma_1)=0,
\end{equation}
\begin{equation}\label{bouss-eq-2}
D^3\alpha_2+3D^2\alpha_3+3D\alpha_4+\gamma_1D\alpha_2-\alpha_1D^2\gamma_1+
2\alpha_2D\gamma_1+\alpha_1D\gamma_2=0,
\end{equation}
respectively.

Note that $\res(Q^{n+m/3})=\alpha_4-\alpha_1D\gamma_1+\alpha_2\gamma_1+\alpha_1\gamma_2$.
From the Gelfand--Dikii equations, it follows that
\begin{eqnarray}
\frac{i\sqrt{3}}{c_{n-1,m}}\frac{\d\gamma_1}{\d t_{n-1,m}} & = & 3D\alpha_1,
\notag \\
\frac{i\sqrt{3}}{c_{n-1,m}}\frac{\d\gamma_2}{\d t_{n-1,m}} & = &
3D^2\alpha_1+3D\alpha_2,  \label{bouss-gd} \\
\frac{i\sqrt{3}}{c_{n,m}}\frac{\d\gamma_1}{\d t_{n,m}} & = &
3D\,\res(Q^{n+m/3}). \notag
\end{eqnarray}

Using (\ref{bouss-eq-1}), (\ref{bouss-eq-2}) we can express $D\alpha_4$ in terms of $\alpha_1$ and $\alpha_2$. Using the first two equations in~(\ref{bouss-gd}) and the string equation, we can express $\alpha_1$ and $\alpha_2$ in terms of derivatives of $F$. If we replace all $\alpha_i$ and
$\gamma_i$ in the third equation in~(\ref{bouss-gd}) with their expressions, we obtain the following:
\begin{multline}\label{formula-bouss}
(3n+m+1)\ll\tau_{n,m}\tau_{0,0}^2\rr =
\ll\tau_{n-1,m}\tau_{0,1}\rr\ll\tau_{0,0}^3\rr+\\
2\ll\tau_{n-1,m}\tau_{0,0}\rr\ll\tau_{0,1}\tau_{0,0}^2\rr+
 2\ll\tau_{n-1,m}\tau_{0,1}\tau_{0,0}\rr\ll\tau_{0,0}^2\rr+\\
\frac{2}{3}\ll\tau_{n-1,m}\tau_{0,1}\tau_{0,0}^3\rr+
 3\ll\tau_{n-1,m}\tau_{0,0}^2\rr\ll\tau_{0,1}\tau_{0,0}\rr.
\end{multline}

This equation allows us to calculate the coefficients of $F$ in the examples below. The similar equation is used in~\cite{w} in the case of $r=2$. It is easy to see that, using the same argument, one can obtain the similar equation for any Gelfand--Dikii hierarchy.




\section{An algorithm to calculate
$\l\tau_{n,m}\prod_{i=1}^{r-1}\tau_{0,i}^{k_{i}}\r_g$}

We fix $r\geq 2$. In this section, we expalain how to calculate intersection numbers from the previous section in some particular cases. We describe a way to calculate $\l\tau_{n,m}\prod_{i=1}^{r-1}\tau_{0,i}^{k_{i}}\r_g$, which looks like follows. First, there is a generalization of our Lemma~\ref{lemma-proof2}. We express intersection numbers $\l\tau_{n,m}\prod_{i=1}^{r-1}\tau_{0,i}^{k_{i}}\r_g$ via integrals against two-pointed ramification cycles (see Section~\ref{tprc} for the definitions of the integrals against two-pointed ramification cycles and Section~\ref{first-step} for the expression). These integrals can be considered as a generalization of the notion of Hurwitz numbers. Then we introduce a recursion relation for these integrals (Section~\ref{recurela}). This relation can be considered as a generalization of a standard recursion relation for Hurwitz numbers, see~\cite{op2}. Thus we
express our integrals via the simplest possible ones, defined only in genera zero and one (Section~\ref{inival}).

The integrals in genus zero equal to the initial intersection numbers $\l\prod_{i=1}^{r-1}\tau_{0,i}^{k_{i}}\r_0$ in genus zero, which we suppose to be known. And the integrals in genus one can be calculated using the topological recursion relation (the last is not completely clear, so we give only the general idea of such calculations, see Section~\ref{genusone}).

For the simple examples for our algorithm, see Section~\ref{simplex}. Another example is computed in the appendix. The proofs of the theorems of this section are collected in the Section~13.

\subsection{Two-pointed ramification cycles}\label{tprc}

We introduce the integrals against two-pointed ramification cycles which we will use in our formulas.

By $V^*_{g,m}(\prod_{i=1}^s\eta_{q_i,a_i})$ denote the subvariety of $\M'_{g,1+s}$ consisting of curves $(C,x_1\dots,x_{1+s},\T)$ such that $-(\sum_{i=1}^sa_i)x_1+\sum_{i=1}^sa_ix_{1+i}$ is the divisor of a meromorphic function. The covering $\M'_{g,1+s}\to\M_{g,1+s}$ is defined here by $m_1=m$, $m_2=q_1,\dots,m_{1+s}=q_s$. All $a_i$ are supposed to be positive integers and also we require $0\leq m,q_1,\dots,q_s\leq r-1$.

We denote by $V_{g,m}(\prod_{i=1}^s\eta_{q_i,a_i})$ the closure of $V^*_{g,m}(\prod_{i=1}^s\eta_{q_i,a_i})$ and denote by $S^n_{g,m}(\prod_{i=1}^s\eta_{q_i,a_i})$ the intersection number
\begin{equation}
\frac{1}{r^g}\int_{V_{g,m}(\prod_{i=1}^s\eta_{q_i,a_i})}
\psi(x_1)^n\cdot c_D(\mathcal{V}).
\end{equation}

Note that $S^{n}_{g,m}(\prod_{i=1}^s\eta_{m_i,a_i})$ is defined if and only if $g\geq 0$; $0\leq m,q_1,\dots,q_s\leq r-1$; $n\geq 0$; $s\geq 1$; and $a_1,\dots,a_s\geq 1$. Moreover, for convenience, we put $S^{n}_{0,m}(\eta_{m_1,1})=0$; and if $g<0$, then we also put $S^{n}_{g,m}(\prod_{i=1}^s\eta_{m_i,a_i})=0$.

Another definition we need is the following one. Consider the moduli space $\oM'_{1,k}\ni(C,x_1,\dots,x_k,\T)$ determined by some labels $m_1,\dots,m_k$. Let $b_1,\dots,b_k$ be nonzero integers such that $\sum_{t=1}^kb_t=0$. By $W(b_1,\dots,b_k)$ denote the closure of the subvariety consisting of smooth curves $(C,x_1,\dots,x_k,\T)$ such that there exists a meromorphic function whose divisor is equal to $\sum_{t=1}^kb_tx_t$.

By $\wS_1(\prod_{t=1}^k\eta_{m_t,b_t})$ denote $(1/r)\int_{W(b_1,\dots,b_k)}c_D(\V)$. We will discuss later how to calculate $\wS_1(\prod_{t=1}^k\eta_{m_t,b_t})$.

\subsection{The first step}\label{first-step}

We state the formula expresssing the intersection numbers $\l\tau_{n,m}\prod_{i=1}^{r-1}\tau_{0,i}^{k_{i}}\r_g$ via the intersection numbers $S^n_{g,m}(\prod_{i=1}^s\eta_{q_i,a_i})$.

\begin{theorem}\label{algorithm1}
If $s\geq 1$, then
\begin{equation}\label{alg1}
\l\tau_{n,m}\prod_{i=1}^{s}\tau_{0,m_i}\r_g=
\sum_{j=0}^g\frac{(-1)^j}{g!}\binom{g}{j}
S^{n-j}_{g,m}(\prod_{i=1}^s\eta_{m_i,1}\cdot\eta_{0,1}^{g-j}).
\end{equation}
\end{theorem}

This theorem can be generalized (one can add a parameter $l$ as it is done in Lemma~\ref{lemma-proof2}).

\subsection{The recursive relation}\label{recurela}

We state our recursive relation for the numbers $S^n_{g,m}(\prod_{i=1}^s\eta_{q_i,a_i})$.

\begin{theorem}\label{algorithm2}
If $n\geq 1$, then
\begin{multline}\label{alg2}
S^{n}_{g,m}(\prod_{i=1}^s\eta_{m_i,a_i})=\\
\sum_{I\subset\{1,\dots,s\}}
\sum_{j=1}^{a(I)}
\sum_{B(j,a(I))}
\left(
\frac{|I|+j-2}{(\sum_{r=1}^sa_r)\cdot (2g+s-1)}
\cdot
\frac{\prod_{r=1}^jb_r}{\aut(b_1,\dots,b_j)}
\cdot
\right.
\\
\sum_{u_1,\dots,u_j=0}^{r-2}
S^{n-1}_{g+1-j,m}(\prod_{i\not\in I}\eta_{m_i,a_i}\prod_{t=1}^j\eta_{u_t,b_t})
\cdot
\l\prod_{t=1}^j\tau_{0,r-2-u_t}\prod_{i\in I}\tau_{0,m_i}\r_0
\\
+
\frac{|I|+j}{(\sum_{r=1}^sa_r)\cdot (2g+s-1)}
\cdot
\frac{\prod_{r=1}^jb_r}{\aut(b_1,\dots,b_j)}
\cdot
\\
\left.
\sum_{u_1,\dots,u_j=0}^{r-2}
S^{n-1}_{g-j,m}(\prod_{i\not\in I}\eta_{m_i,a_i}\prod_{t=1}^j\eta_{u_t,b_t})
\cdot
\wS_1(\prod_{i\not\in I}\eta_{m_i,-a_i}\prod_{t=1}^j\eta_{r-2-u_t,b_t})
\right).
\end{multline}
Here the first sum is taken over all subsets $I$ of $\{1,\dots,s\}$. Then, $a(I)$ is defined to be $\sum_{k\in I}a_k$. The third sum is taken over all possible partitions $B=(b_1,\dots,b_j)$ of $a(I)$ of length $j$; $\sum_{k=1}^jb_k=a(I)$, $b_1\geq\dots\geq b_j$, and all $b_k$ are positive integers.
\end{theorem}

The examples of applying this recursive relation can be found in Section 12.5 and in the appendix.

\subsection{Initial values}\label{inival}

The initial values look as in the following theorem.

\begin{theorem}\label{algorithm3}
We have
\begin{eqnarray}
S^{0}_{g,m}(\prod_{i=1}^s\eta_{m_i,a_i}) & = & 0,\quad if\ g>1; \\
S^0_{1,m}(\prod_{i=1}^s\eta_{m_i,a_i}) & = &
\wS_1(\eta_{m,-\sum_{i=1}^sa_i}\cdot\prod_{i=1}^s\eta_{m_i,a_i}); \\
S^{n}_{0,m}(\prod_{i=1}^s\eta_{m_i,a_i}) & = & \l\tau_{n,m}
\prod_{i=1}^s\tau_{0,m_i}\r_0.
\end{eqnarray}
\end{theorem}

So, applying (\ref{alg1}) and then several times (\ref{alg2}), we express the intersection number
$\l\tau_{n,m}\prod_{i=1}^{r-1}\tau_{0,i}^{k_{i}}\r_g$ via the intersection numbers
$\l\prod_{i=1}^{r-1}\tau_{0,i}^{p_{i}}\r_0$ and $\wS_1(\prod_{t=1}^k\eta_{m_t,b_t})$. Later we will explain how to express the intersection numbers $\wS_1(\prod_{t=1}^k\eta_{m_t,b_t})$ via $\l\prod_{i=1}^{r-1}\tau_{0,i}^{p_{i}}\r_0$.

\subsection{Simple examples}\label{simplex}

Here, we use Theorems~\ref{algorithm1}, \ref{algorithm2}, and~\ref{algorithm3} to compute the Mumfor--Morita--Miller intersection numbers in some special cases.

\subsubsection{$\l\tau_{1,0}\r_1$}
We consider the case $r=4$. From the topological recursion relation (see~\cite{g}), we know that $\l\tau_{1,0}\r_1=1/8$. We can prove this independently.

Note that $\l\tau_{1,0}\r_1=\l\tau_{2,0}\tau_{0,0}\r_1$. Equation~(\ref{alg1}) implies that
$\l\tau_{2,0}\tau_{0,0}\r_1=S^{2}_{1,0}(\eta_{0,1}^2)$.

We know from~\cite{w1} that, in the case $r=4$, $\l\tau_{0,0}^2\tau_{0,2}\r_0=\l\tau_{0,0}\tau_{0,1}^2\r_0=1$, $\l\tau_{0,1}^2\tau_{0,2}^2\r_0=1/4$, $\l\tau_{0,2}^5\r_0=1/8$, and all other $\l\prod_{i=1}^{r-1}\tau_{0,i}^{p_{i}}\r_0$ are equal to zero. Another fact we need here is that $\wS_1(\prod_{t=1}^k\eta_{m_t,b_t})$ is equal to zero if one $m_i$ is equal to zero.

Then, from~(\ref{alg2}) and formulas for initial values, it follows that
\begin{eqnarray*}
S^{2}_{1,0}(\eta_{0,1}^2) & = &
\frac{1}{3}S^1_{1,0}(\eta_{0,2})\l\tau_{0,0}^2\tau_{0,2}\r_0;\\
S^1_{1,0}(\eta_{0,2}) & = & \frac{1}{4}S^0_{0,0}(\eta_{0,1}\eta_{2,1})
\l\tau_{0,0}^2\tau_{0,2}\r_0+
\frac{1}{8}S^0_{0,0}(\eta_{1,1}^2)
\l\tau_{0,0}\tau_{0,1}^2\r_0;\\
S^0_{0,0}(\eta_{0,1}\eta_{2,1}) & = & \l\tau_{0,0}^2\tau_{0,2}\r_0;\\
S^0_{0,0}(\eta_{1,1}^2) & = & \l\tau_{0,0}\tau_{0,1}^2\r_0.
\end{eqnarray*}

Thus we obtain that $\l\tau_{1,0}\r_1=1/8$ as it has to be.

\subsubsection{$\l\tau_{1,1}\tau_{0,1}^3\tau_{0,0}\r_0$}

We calculate $\l\tau_{1,1}\tau_{0,1}^3\tau_{0,0}\r_0$ in the case $r=3$. It follows from the
string equation that $\l\tau_{1,1}\tau_{0,1}^3\tau_{0,0}\r_0=\l\tau_{0,1}^4\r_0=1/3$, but we want to calculate this using our algorithm.

Recall that in the case $r=3$, $\l\tau_{0,0}^2\tau_{0,1}\r_0=1$, $\l\tau_{0,1}^4\r_0=1/3$, and all other $\l\prod_{i=1}^{r-1}\tau_{0,i}^{p_{i}}\r_0$ are equal to zero.

From our algorithm, we have the following
\begin{eqnarray*}
\l\tau_{1,1}\tau_{0,1}^3\tau_{0,0}\r_0 & = & S^1_{0,1}(\eta_{1,1}^3\eta_{0,1})
\\
S^1_{0,1}(\eta_{1,1}^3\eta_{0,1}) &
= & \frac{1}{2}S^0_{0,1}(\eta_{1,1}^2\eta_{1,2})
\l\tau_{0,0}^2\tau_{0,1}\r_0+\frac{1}{2}S^0_{0,1}(\eta_{0,3}\eta_{0,1})
\l\tau_{0,1}^4\r_0 \\
S^0_{0,1}(\eta_{1,1}^2\eta_{1,2}) &= &\l\tau_{0,1}^4\r_0 \\
S^0_{0,1}(\eta_{0,3}\eta_{0,1}) & = & \l\tau_{0,0}^2\tau_{0,1}\r_0
\end{eqnarray*}

Thus we obtain that $\l\tau_{1,1}\tau_{0,1}^3\tau_{0,0}\r_0=1/3$.

\subsection{Calculation of $\wS_1(\prod_{t=1}^k\eta_{m_t,b_t})$}\label{genusone}

Here, we explain only a general idea how to calculate any number $\wS_1(\prod_{t=1}^k\eta_{m_t,b_t})$. Then we give the concrete calculations in the case $r=4$. It is the first case where such numbers are not identically zero.

\subsubsection{Our idea}

There is a topological recursion relation expressing intersection numbers in genus one via intersection numbers in genus zero (see, e.~g., \cite{jkv}). We need the special case of this relation, which looks as follows:
\begin{equation}
\l\tau_{n+1,m}\prod_{i=1}^s\tau_{0,m_i}\r_1=
\frac{1}{24}\sum_{l=0}^{r-2}
\l\tau_{n,m}\tau_{0,l}\tau_{0,r-2-l}\prod_{i=1}^s\tau_{0,m_i}\r_0.
\end{equation}

Using this relation we can calculate all intersection numbers $\l\tau_{n+1,m}\prod_{i=1}^s\tau_{0,m_i}\r_1$ (of course, we always suppose that we know all intersection numbers in genus zero). Then we can try to calculate the same intersection numbers using our algorithm based on Theorems~\ref{algorithm1}, \ref{algorithm2}, and \ref{algorithm3}. This gives us a number of linear equations for the numbers $\wS_1(\prod_{t=1}^k\eta_{m_t,b_t})$. Then we just have to solve these equations.

Of course, it is not obvious that this way allows to calculate the numbers $\wS_1(\prod_{t=1}^k\eta_{m_t,b_t})$ in the case of arbitrary $r$. Nevertheless, we are almost sure that this works for any $r$.

\subsubsection{Example}

Consider the case $r=4$. It is easy to see that $\wS_1(\prod_{t=1}^k\eta_{m_t,b_t})\not=0$ if and only if $k=2$, $m_1=m_2=2$, and $b_1=-b_2$. We calculate $\wS_1(\eta_{2,2}\eta_{2,-2})$ and $\wS_1(\eta_{2,3}\eta_{2,-3})$.

Note that
\begin{equation}
\frac{1}{3\cdot 2^5}=
\l\tau_{1,2}\tau_{0,2}\r_1=
S^1_{1,2}(\eta_{2,1}\eta_{0,1})=
\frac{1}{3}
\wS_1(\eta_{2,2}\eta_{2,-2}).
\end{equation}
Therefore, $\wS_1(\eta_{2,2}\eta_{2,-2})=1/2^5$.

Then
\begin{eqnarray}
\l\tau_{1,2}\tau_{0,2}\r_1 & = &
S^1_{1,2}(\eta_{2,2}\eta_{0,1})-\wS_1(\eta_{2,2}\eta_{2,-2});\\
S^1_{1,2}(\eta_{2,2}\eta_{0,1}) & = &
\frac{1}{3}\wS_1(\eta_{2,3}\eta_{2,-3})+
\frac{4}{9}\wS_1(\eta_{2,2}\eta_{2,-2}).
\end{eqnarray}
Therefore, $\wS_1(\eta_{2,3}\eta_{2,-3})=1/12$.




\section{Proofs of Theorems~\ref{algorithm1}, \ref{algorithm2}, and~\ref{algorithm3}}

\subsection{Initial values}\label{initial}

\begin{proof}[Proof of Theorem~\ref{algorithm3}]
Consider the intersection number $S=S^{0}_{g,m}(\prod_{i=1}^s\eta_{m_i,a_i})$, where $g>1$. We want to prove that $S=0$. Recall that $S$ is defined as follows:
\begin{equation}
S=\frac{1}{r^g}\int_{V_{g,m}(\prod_{i=1}^s\eta_{q_i,a_i})}c_D(\mathcal{V})
\end{equation}
Note that $\dim V_{g,m}(\prod_{i=1}^s\eta_{q_i,a_i})=2g+s-2$ and $D<g+s$. If $g\geq 2$, then $2g+s-2\geq g+s$, and we obtain $S=0$.

Since $V_{1,m}(\prod_{i=1}^s\eta_{q_i,a_i})=W(-\sum_{i=1}^sa_i,a_1,\dots,a_s)$, where $m_1=m$, $m_2=q_1,\dots,m_{s+1}=q_s$, if follows that
\begin{equation}
S^0_{1,m}(\prod_{i=1}^s\eta_{m_i,a_i})=
\wS_1(\eta_{m,-\sum_{i=1}^sa_i}\cdot\prod_{i=1}^s\eta_{m_i,a_i}).
\end{equation}

The equality $S^{n}_{0,m}(\prod_{i=1}^s\eta_{m_i,a_i})=\l\tau_{n,m}\prod_{i=1}^s\tau_{0,m_i}\r_0$ is obtained just from the fact that, in the case of genus zero, $V_{0,m}(\prod_{i=1}^s\eta_{q_i,a_i})$ is equal to the appropriate space $\oM'_{0,1+s}$.
\end{proof}

\subsection{First step of the algorithm}

\begin{proof}[Proof of Theorem~\ref{algorithm1}]
Theorem~\ref{algorithm1} is proved by the very same argument as Lemma~\ref{lemma-proof2}. The only difference is the following one. We have to formulate and prove an analogue of Lemma~\ref{S-expression-via-difference}. Proving this analogue and using the lemma of E.~Ionel, we represent the $\psi$-class as a sum of divisors. Then we choose only those divisors, where $\psi(y)^k\cdot c_D(\V)$ is not zero. Here we have to use one more additional argument:
$\l\tau_{0,0}\cdot\prod_{i=1}^l\tau_{0,m_i}\r_0\not=0$ if and only if $l=2$ and $m_1+m_2=r-2$. Then this intersection number is equal to $1$. All other steps of the proof are just the same. \end{proof}

\subsection{Recursion relation}

\begin{proof}[Proof of Theorem~\ref{algorithm2}]
This is also proved in the same way as Lemma~\ref{lemma-proof2}. We use the lemma of E.~Ionel. In the target moduli space of the $\hll$ mapping, we express the $\psi$-class as the divisor whose generic point is represented by a two-component curve such that the point corresponding to $x_1$
lies on the first component and the point corresponding to $x_2,\dots,x_s$ with a fixed critical value lies on the other component.

Then we express $\psi(x_1)$ as a sum of some divisors in $\hH$. The mapping $\sigma\circ st$ ($\sigma$ is the projection forgetting all marked point except for $x_1,\dots,x_{1+s}$) takes each divisor to a subvariety of $\pi(V_{g,m}(\prod_{i=1}^s\eta_{q_i,a_i}))$ ($\pi$ is the projection $\oM'_{g,1+s}\to\oM_{g,1+s}$). We only need subvarieties of codimension one. This condition means that all critical points of the corresponding functions are lying exactly on two components of a curve in $st(\hH)$. In other words, this means that a curve in $\sigma(st(\hH))$ consists of two components.

Consider such irreducible divisor in $\pi(V_{g,m}(\prod_{i=1}^s\eta_{q_i,a_i}))=\sigma(st(\hH))$. Two components of a curve representing a generic point of this divisor can intersect at $j$ points. One component contains points $x_{i+1}$, $i\in I\subset\{1,\dots,s\}$, and $j$ points of intersection. The other component contains points $x_1$ and $x_{i+1}$, $i\not\in I$, and also $j$ points of intersection. The first component determines a two-pointed ramification cycle, where the divisor is $\sum_{i\in I}a_ix_{1+i}-\sum_{t=1}^jb_t*_t$ (by $*_t$ denote the points of intersection). The other component also determines a two-pointed ramification cycle, where the
divisor is $-(\sum_{t=1}^sa_t)x_1+\sum_{i\not\in I}a_ix_{1+i}+\sum_{t=1}^jb_t*_t$.

Let the first component have genus $g_1$ and the second one have genus $g_2$. We have $g_1+g_2+j-1=g$. We must consider the preimage of this divisor under the mapping $\pi$ and then integrate against it the class $\psi(x_1)^{n-1}\cdot c_D(\V)$. Note that when $c_D(\V)$ does not vanish, it factorizes to $c_{D_1}(\V_1)\cdot c_{D_2}(\V_2)$. Then we have to integrate $c_{D_1}(\V_1)$ over the two-pointed ramification cycle determined by the first component, and
$\psi(x_1)^{n-1}\cdot c_{D_2}(\V_2)$ over the two-pointed ramification cycle determined by the second component.

In the first case, we see that dimensional conditions (as in Section~\ref{initial}) imply that the integral does not vanish if and only if $g=1$ or $g=0$. Thus we obtain the second multipliers in the formula. The integral corresponding to the second component obviously gives us the first
multipliers in the formula.

Now we only have to explain the coefficients in the formula. The coefficient $1/r^g$ appearing in the definition of $S^{n}_{g,m}(\prod_{i=1}^s\eta_{m_i,a_i})$ behaves properly since $\pi$ is a
ramified covering with correspoding multiplicities. Then $\prod_{r=1}^jb_r$ is the multiplicity of $\hll$ along the corresponding divisor in $\hH$; the coefficient $\aut(b_1,\dots,b_j)$ appears
since we have to mark the points of intersection of two components; $(\sum_{r=1}^sa_r)$ comes from the lemma of E.~Ionel; and $(|I|+j-2)/(2g+s-1)$ in the case of genus $0$ (or $(|I|+j)/(2g+s-1)$ in the case of genus onr) is the fraction of multiplicities of $\sigma\circ st$ over the divisor and over the initial subvariety $\pi(V_{g,m}(\prod_{i=1}^s\eta_{q_i,a_i}))$.

Thus we obtain the required formula.
\end{proof}




\appendix

\section{Calculation of $\l\tau_{6,1}\r_3$ in the case $r=3$}

In this section, we calculate the intersection number $\l\tau_{6,1}\r_3$ in the case $r=3$ using our algorithm (Theorems~\ref{algorithm1}, \ref{algorithm2}, and~\ref{algorithm3}). Then we calculate the corresponding coefficient of the string solution of the Boussinesq hierarchy using the relation explained in Section~\ref{relation}. The results will appear to be the same. Thus we check a very special case of the Witten's conjecture.

\subsection{First step of the algorithm}
The first step is the following:
\begin{equation*}
\l\tau_{6,1}\r_3=\l\tau_{7,1}\tau_{0,0}\r_3=
\frac{1}{6}S^7_{3,1}(\eta_{0,1}^4)-\frac{1}{2}S^6_{3,1}(\eta_{0,1}^3)
+\frac{1}{2}S^5_{3,1}(\eta_{0,1}^2).
\end{equation*}

In the next three subsections, we calculate separately the summands of this expression. We recall once again that, in the case of $r=3$, $\l\tau_{0,0}^2\tau_{0,1}\r_0=1$, $\l\tau_{0,1}^4\r_0=1/3$, and all other $\l\prod_{i=1}^{r-1}\tau_{0,i}^{p_{i}}\r_0$ are equal to zero.

For convenience, we shall denote $\l\tau_{0,0}^2\tau_{0,1}\r_0$ by $\l 1\r$ and $\l\tau_{0,1}^4\r_0$ by $\l 1/3\r$.

\subsection{Calculations in degree $4$}

Using (\ref{alg2}), we get
\begin{equation*}
S^7_{3,1}(\eta_{0,1}^4)=\frac{1}{3}S^6_{3,1}(\eta_{0,1}^2\eta_{0,2})
\l 1\r.
\end{equation*}

Then
\begin{eqnarray*}
S^6_{3,1}(\eta_{0,1}^2\eta_{0,2}) &
= & \frac{3}{16}S^5_{3,1}(\eta_{0,1}\eta_{0,3})\l 1\r+
+\frac{1}{16}S^5_{3,1}(\eta_{0,2}^2)\l 1\r+
\frac{1}{32}S^5_{2,1}(\eta_{0,1}^3\eta_{1,1})\l 1\r;\\
S^5_{3,1}(\eta_{0,1}\eta_{0,3}) & = &
\frac{1}{7}S^4_{3,1}(\eta_{0,4})\l 1\r+
\frac{1}{14}S^4_{2,1}(\eta_{0,1}^2\eta_{1,2})\l 1\r+
\frac{1}{14}S^4_{2,1}(\eta_{0,1}\eta_{0,2}\eta_{1,1})\l 1\r; \\
S^5_{3,1}(\eta_{0,2}^2) & = &
\frac{1}{7}S^4_{3,1}(\eta_{0,4})\l 1\r+
\frac{1}{14}S^4_{2,1}(\eta_{0,1}\eta_{0,2}\eta_{1,1})\l 1\r; \\
S^5_{2,1}(\eta_{0,1}^3\eta_{1,1}) & = &
\frac{3}{14}S^4_{2,1}(\eta_{0,1}^2\eta_{1,2})\l 1\r+
\frac{3}{14}S^4_{2,1}(\eta_{0,1}\eta_{0,2}\eta_{1,1})\l 1\r.
\end{eqnarray*}

Therefore,
\begin{equation*}
S^7_{3,1}(\eta_{0,1}^4)=
\frac{1}{84}S^4_{3,1}(\eta_{0,4})+
\frac{3}{448}S^4_{2,1}(\eta_{0,1}^2\eta_{1,2})+
\frac{11}{1344}S^4_{2,1}(\eta_{0,1}\eta_{0,2}\eta_{1,1}).
\end{equation*}

Then
\begin{eqnarray*}
S^4_{3,1}(\eta_{0,4}) & = &
\frac{1}{8}S^3_{2,1}(\eta_{0,1}\eta_{1,3})\l 1\r+
\frac{1}{8}S^3_{2,1}(\eta_{0,3}\eta_{1,1})\l 1\r+
\frac{1}{6}S^3_{2,1}(\eta_{0,2}\eta_{1,2})\l 1\r;\\
S^4_{2,1}(\eta_{0,1}^2\eta_{1,2}) & = &
\frac{1}{4}S^3_{2,1}(\eta_{0,1}\eta_{1,3})\l 1\r+
\frac{1}{12}S^3_{2,1}(\eta_{0,2}\eta_{1,2})\l 1\r+
\frac{1}{48}S^3_{1,1}(\eta_{0,1}^2\eta_{1,1}^2)\l 1\r;\\
S^4_{2,1}(\eta_{0,1}\eta_{0,2}\eta_{1,1}) & = &
\frac{1}{8}S^3_{2,1}(\eta_{0,1}\eta_{1,3})\l 1\r+
\frac{1}{12}S^3_{2,1}(\eta_{0,2}\eta_{1,2})\l 1\r+ \\
& &
\frac{1}{8}S^3_{2,1}(\eta_{0,3}\eta_{1,1})\l 1\r+
\frac{1}{24}S^3_{1,1}(\eta_{0,1}^2\eta_{1,1}^2)\l 1\r.
\end{eqnarray*}

Therefore,
\begin{multline*}
S^7_{3,1}(\eta_{0,1}^4)=
\frac{15}{2^9\cdot 7}
S^3_{2,1}(\eta_{0,1}\eta_{1,3})+
\frac{13}{2^6\cdot 3^2\cdot 7}
S^3_{2,1}(\eta_{0,2}\eta_{1,2})+\\
\frac{9}{2^9\cdot 7}
S^3_{2,1}(\eta_{0,3}\eta_{1,1})+
\frac{31}{2^{10}\cdot 3^2\cdot 7}
S^3_{1,1}(\eta_{0,1}^2\eta_{1,1}^2).
\end{multline*}

Note that
\begin{eqnarray*}
S^3_{2,1}(\eta_{0,1}\eta_{1,3}) & = &
\frac{1}{5} S^2_{2,1}(\eta_{1,4}) \l 1 \r +
\frac{1}{10} S^2_{1,1}(\eta_{0,1}\eta_{1,1}\eta_{1,2}) \l 1 \r +
\frac{1}{60} S^2_{0,1}(\eta_{0,0}^4) \l \frac{1}{3} \r; \\
S^3_{2,1}(\eta_{0,2}\eta_{1,2}) & = &
\frac{1}{5} S^2_{2,1}(\eta_{1,4}) \l 1 \r +
\frac{1}{20} S^2_{1,1}(\eta_{0,1}\eta_{1,1}\eta_{1,2}) \l 1 \r +
\frac{1}{40} S^2_{1,1}(\eta_{0,2}\eta_{1,1}^2) \l 1 \r; \\
S^3_{2,1}(\eta_{0,3}\eta_{1,1}) & = &
\frac{1}{5} S^2_{2,1}(\eta_{1,4}) \l 1 \r +
\frac{1}{10} S^2_{1,1}(\eta_{0,1}\eta_{1,1}\eta_{1,2}) \l 1 \r +
\frac{1}{10} S^2_{1,1}(\eta_{0,2}\eta_{1,1}^2) \l 1 \r; \\
S^3_{1,1}(\eta_{0,1}^2\eta_{1,1}^2) & = &
\frac{2}{5} S^2_{1,1}(\eta_{0,1}\eta_{1,1}\eta_{1,2}) \l 1 \r +
\frac{1}{10} S^2_{1,1}(\eta_{0,2}\eta_{1,1}^2) \l 1 \r +
\frac{1}{20} S^2_{0,1}(\eta_{0,0}^4) \l \frac{1}{3} \r.
\end{eqnarray*}

Then
\begin{eqnarray*}
S^2_{2,1}(\eta_{1,4}) & = &
\frac{3}{16} S^1_{1,1}(\eta_{1,3}\eta_{1,1}) \l 1 \r +
\frac{1}{8} S^1_{1,1}(\eta_{1,2}^2) \l 1 \r 
+\frac{1}{8} S^1_{0,1}(\eta_{0,1}^2\eta_{0,2}) \l \frac{1}{3} \r; \\
S^2_{1,1}(\eta_{0,1}\eta_{1,1}\eta_{1,2}) & = &
\frac{3}{16} S^1_{1,1}(\eta_{1,3}\eta_{1,1}) \l 1 \r +
\frac{1}{8} S^1_{1,1}(\eta_{1,2}^2) \l 1 \r + \\
& &
\frac{1}{4} S^1_{0,1}(\eta_{0,1}^2\eta_{0,2}) \l \frac{1}{3} \r +
\frac{1}{32} S^1_{0,1}(\eta_{0,1}\eta_{1,1}^3) \l 1 \r; \\
S^2_{1,1}(\eta_{0,2}\eta_{1,1}^2) & = &
\frac{3}{8} S^1_{1,1}(\eta_{1,3}\eta_{1,1}) \l 1 \r +
\frac{1}{16} S^1_{0,1}(\eta_{0,1}\eta_{1,1}^3) \l 1 \r +
\frac{1}{16} S^1_{0,1}(\eta_{0,1}^2\eta_{0,2}) \l \frac{1}{3} \r.
\end{eqnarray*}

Since
\begin{eqnarray*}
S^1_{1,1}(\eta_{1,3}\eta_{1,1}) & = &
\frac{1}{6} S^0_{0,1}(\eta_{1,1}^2\eta_{1,2}) \l 1 \r +
\frac{1}{2} S^0_{0,1}(\eta_{0,1}\eta_{0,3}) \l \frac{1}{3} \r+
\frac{1}{3} S^0_{0,1}(\eta_{0,2}^2) \l \frac{1}{3} \r \\
& = & \frac{1}{3};\\
S^1_{1,1}(\eta_{1,2}^2) & = &
\frac{1}{12} S^0_{0,1}(\eta_{1,1}^2\eta_{1,2}) \l 1 \r +
\frac{1}{2} S^0_{0,1}(\eta_{0,1}\eta_{0,3}) \l \frac{1}{3} \r+
\frac{1}{3} S^0_{0,1}(\eta_{0,2}^2) \l \frac{1}{3} \r  \\
& = & \frac{11}{2^2\cdot 3^2},
\end{eqnarray*}
it follows that
\begin{eqnarray*}
S^2_{2,1}(\eta_{1,4}) & = & \frac{41}{2^5\cdot 3^2};\\
S^2_{1,1}(\eta_{0,1}\eta_{1,1}\eta_{1,2}) & = & \frac{7}{2^2\cdot 3^2};\\
S^2_{1,1}(\eta_{0,2}\eta_{1,1}^2) & = & \frac{1}{2\cdot 3};
\end{eqnarray*}
and, therefore,
\begin{eqnarray*}
S^3_{2,1}(\eta_{0,1}\eta_{1,3}) & = & \frac{77}{2^5\cdot 3^2\cdot 5}\\
S^3_{2,1}(\eta_{0,2}\eta_{1,2}) & = & \frac{61}{2^5\cdot 3^2\cdot 5}\\
S^3_{2,1}(\eta_{0,3}\eta_{1,1}) & = & \frac{31}{2^5\cdot 3\cdot 5}\\
S^3_{1,1}(\eta_{0,1}^2\eta_{1,1}^2) & = & \frac{1}{3^2}.
\end{eqnarray*}

Thus we have
\begin{equation*}
S^7_{3,1}(\eta_{0,1}^4)=\frac{209}{2^7\cdot 3^4\cdot 5\cdot 7}.
\end{equation*}

\subsection{Calculations in degree $3$}

We have
\begin{eqnarray*}
S^6_{3,1}(\eta_{0,1}^3) & = & \frac{1}{4} S^5_{3,1}(\eta_{0,1}\eta_{0,2}) \l 1 \r; \\
S^5_{3,1}(\eta_{0,1}\eta_{0,2}) & = &
\frac{1}{7} S^4_{3,1}(\eta_{0,3}) \l 1 \r +
\frac{1}{21} S^4_{2,1}(\eta_{0,1}^2\eta_{1,1}) \l 1 \r; \\
S^4_{3,1}(\eta_{0,3}) & = &
\frac{1}{9} S^3_{2,1}(\eta_{0,1}\eta_{1,2}) \l 1 \r +
\frac{1}{9} S^3_{2,1}(\eta_{0,2}\eta_{1,1}) \l 1 \r; \\
S^4_{2,1}(\eta_{0,1}^2\eta_{1,1}) & = &
\frac{2}{9} S^3_{2,1}(\eta_{0,1}\eta_{1,2}) \l 1 \r +
\frac{1}{9} S^3_{2,1}(\eta_{0,2}\eta_{1,1}) \l 1 \r; \\
S^3_{2,1}(\eta_{0,1}\eta_{1,2}) & = &
\frac{1}{5} S^2_{2,1}(\eta_{1,3}) \l 1 \r +
\frac{1}{30} S^2_{1,1}(\eta_{0,1}\eta_{1,1}^2) \l 1 \r; \\
S^3_{2,1}(\eta_{0,2}\eta_{1,1}) & = &
\frac{1}{5} S^2_{2,1}(\eta_{1,3}) \l 1 \r +
\frac{1}{15} S^2_{1,1}(\eta_{0,1}\eta_{1,1}^2) \l 1 \r; \\
S^2_{2,1}(\eta_{1,3}) & = &
\frac{1}{6} S^1_{1,1}(\eta_{1,1}\eta_{1,2}) \l 1 \r +
\frac{1}{36} S^1_{0,1}(\eta_{0,1}^3) \l \frac{1}{3} \r; \\
S^2_{1,1}(\eta_{0,1}\eta_{1,1}^2) & = &
\frac{1}{3} S^1_{1,1}(\eta_{1,1}\eta_{1,2}) \l 1 \r +
\frac{1}{12} S^1_{0,1}(\eta_{0,1}^3) \l \frac{1}{3} \r; \\
S^1_{1,1}(\eta_{1,1}\eta_{1,2}) & = &
\frac{4}{9} S^0_{0,1}(\eta_{0,1}\eta_{0,2}) \l \frac{1}{3} \r +
\frac{1}{18} S^0_{0,1}(\eta_{1,1}^3) \l 1 \r =\frac{1}{2\cdot 3}.
\end{eqnarray*}

Thus we have
\begin{eqnarray*}
S^2_{1,1}(\eta_{0,1}\eta_{1,1}^2) & = & \frac{1}{2^2\cdot 3};\\
S^2_{2,1}(\eta_{1,3}) & = & \frac{1}{3^3};\\
S^3_{2,1}(\eta_{0,2}\eta_{1,1}) & = & \frac{7}{2^2\cdot 3^3\cdot 5};\\
S^3_{2,1}(\eta_{0,1}\eta_{1,2}) & = & \frac{11}{2^3\cdot 3^3\cdot 5};\\
S^4_{2,1}(\eta_{0,1}^2\eta_{1,1}) & = & \frac{1}{2\cdot 3^3\cdot 5};\\
S^4_{3,1}(\eta_{0,3}) & = & \frac{5}{2^3\cdot 3^5};\\
S^5_{3,1}(\eta_{0,1}\eta_{0,2}) & = & \frac{37}{2^3\cdot 3^5\cdot 5\cdot 7}
\end{eqnarray*}
and, therefore,
\begin{equation*}
S^6_{3,1}(\eta_{0,1}^3)=\frac{37}{2^5\cdot 3^5\cdot 5\cdot 7}.
\end{equation*}

\subsection{Calculations in degree $2$}

We have
\begin{eqnarray*}
S^5_{3,1}(\eta_{0,1}^2)=
\frac{1}{7}S^4_{3,1}(\eta_{0,2})\l 1\r=
\frac{1}{7}\cdot\frac{1}{12}S^3_{2,1}(\eta_{0,1}\eta_{1,1})\l 1\r=\\
\frac{1}{2^2\cdot 3\cdot 7}\cdot\frac{1}{5}
S^2_{2,1}(\eta_{1,2})\l 1\r=
\frac{1}{2^2\cdot 3\cdot 5\cdot 7}\cdot
\frac{1}{16}S^1_{1,1}(\eta_{1,1}^2)\l 1\r=\\
\frac{1}{2^6\cdot 3\cdot 5\cdot 7}\cdot
\frac{1}{6}S^0_{0,1}(\eta_{0,1}^2)\l \frac{1}{3}\r=
\frac{1}{2^7\cdot 3^3\cdot 5\cdot 7}.
\end{eqnarray*}

\subsection{Summary}

We have
\begin{equation*}
\l\tau_{6,1}\r_3=\frac{209}{2^8\cdot 3^5\cdot 5\cdot 7}
-\frac{37}{2^6\cdot 3^5\cdot 5\cdot 7}
+\frac{1}{2^8\cdot 3^3\cdot 5\cdot 7}
=\frac{1}{2^7\cdot 3^5}.
\end{equation*}

We calculate the corresponding coefficient of the string solution of the
Boussinesq hierarchy. It follows from~(\ref{formula-bouss}) that
\begin{eqnarray*}
26\l\tau_{8,1}\tau_{0,0}^2\r_3 & = & 2\l\tau_{7,1}\tau_{0,0}\r_3
\l\tau_{0,1}\tau_{0,0}^2\r_0+
\frac{2}{3}\l\tau_{7,1}\tau_{0,1}\tau_{0,0}^3\r_2;\\
20\l\tau_{6,1}\tau_{0,1}\tau_{0,0}^2\r_2 & = &
4\l\tau_{5,1}\tau_{0,1}\tau_{0,0}\r_2
\l\tau_{0,1}\tau_{0,0}^2\r_0+
\frac{2}{3}
\l\tau_{5,1}\tau_{0,1}^2\tau_{0,0}^3\r_1;\\
14\l\tau_{4,1}\tau_{0,1}^2\tau_{0,0}^2\r_1 & = &
6\l\tau_{3,1}\tau_{0,1}^2\tau_{0,0}\r_1
\l\tau_{0,1}\tau_{0,0}^2\r_0+
\frac{2}{3}
\l\tau_{4,1}\tau_{0,1}^3\tau_{0,0}^3\r_0.
\end{eqnarray*}

Therefore, $\l\tau_{6,1}\r_3=1/(2^7\cdot 3^5)$. Thus we have checked Witten's conjecture for $\l\tau_{6,1}\r_3$ in the case of Boussinesq
hierarchy.

\end{document}